\documentstyle[amssymb,amsfonts]{amsart}

\def\R{{{\mathbf R}}}

\def\E{{{\mathbf E}}}

\def\P{{{\mathbf P}}}
\def\I{{{\mathbf I}}}
\def\N{{{\mathbf N}}}
\def\Z{{{\mathbf Z}}}

\def\ideal{{\mathbf{i}}}
\def\height{{h}}
\def\B{{\mathcal{B}}}

\def\reg{{\operatorname{reg}}}

\def\max{{\operatorname{max}}}

\def\eps{\varepsilon}
\newenvironment{proof}{\noindent {\bf Proof} }{\endprf\par}
\def \endprf{\hfill  {\vrule height6pt width6pt depth0pt}\medskip}

\theoremstyle{plain}
  \newtheorem{theorem}[subsection]{Theorem}

  \newtheorem{proposition}[subsection]{Proposition}
  
  \newtheorem{lemma}[subsection]{Lemma}
  \newtheorem{corollary}[subsection]{Corollary}

\theoremstyle{remark}
  \newtheorem{remark}[subsection]{Remark}
  \newtheorem{remarks}[subsection]{Remarks}
  \newtheorem{example}[subsection]{Example}
  \newtheorem{examples}[subsection]{Examples}

\theoremstyle{definition}
  \newtheorem{definition}[subsection]{Definition}

\include{psfig}

\begin{document}

\title[Correspondence between graph theory and probability theory]{A correspondence principle between (hyper)graph theory and probability theory, and the (hyper)graph removal lemma}

\author{Terence Tao}
\address{Department of Mathematics, UCLA, Los Angeles CA 90095-1555}
\email{tao@@math.ucla.edu}

\begin{abstract}  We introduce a correspondence principle (analogous to the Furstenberg correspondence principle) that allows one to extract an infinite random graph or hypergraph from a sequence of increasingly large deterministic graphs or hypergraphs.  As an application we present a new (infinitary) proof of the hypergraph removal lemma of Nagle-Schacht-R\"odl-Skokan and Gowers, which does not require the hypergraph regularity lemma and requires significantly less computation.  This in turn gives new proofs of several corollaries of the hypergraph removal lemma, such as Szemer\'edi's theorem on arithmetic progressions. 
\end{abstract}

\maketitle

\section{Introduction}

It is an interesting phenomenon in mathematics that certain types of problems can be treated both by finitary means (e.g. using combinatorial analysis of finite sets), and by infinitary means (e.g. using constructions involving the axiom of choice), thus giving parallel but 
distinct ways to prove a single result.  One particularly striking example of this is \emph{Szemer\'edi's theorem} (see Theorem \ref{sz-quant}) on arithmetic progressions.  This difficult and important theorem now 
has several proofs, both finitary and infinitary, using fields of mathematics as diverse as Fourier analysis, ergodic theory, graph theory, 
hypergraph theory, and elementary combinatorics; the finitary and infinitary arguments are connected by the beautiful \emph{Furstenberg correspondence principle} (see Section \ref{motivation}).  These proofs have different strengths and weaknesses; generally speaking, the infinitary proofs 
are cleaner, shorter, and more elegant, but require
significantly more machinery, whereas the finitary proofs are more elementary and provide more quantitative results, but tend to be messier
and longer in nature.  One particularly visible difference is that finitary proofs often require a number of small parameters (such as $\eps,\delta$)
or large parameters (such as $N, M$), whereas in the infinitary analogues of these proofs, the small parameters often have become zero and 
the large parameters have become infinite, which can lead to cleaner (but more subtle) arguments.

Some progress has been made in reconciling the finitary and infinitary approaches\footnote{From a proof-theoretical perspective, one can use quantifier-elimination methods (such as Herbrand's theorem) to automatically convert a large class of infinitary arguments to finitary ones; this was for instance carried out for the Furstenberg-Weiss infinitary proof of van der Waerden's theorem via topological dynamics, see \cite{girard}.  
However such methods do not seem to shed much light on the connection between the infinitary proofs and the existing finitary proofs in the literature.}, 
as it has been increasingly realized that ideas and methods
from the infinitary world can be transferred to the finitary world, and vice versa; see for instance \cite{tao:ergodic} for a finitary version of
the infinitary ergodic approach to Szemer\'edi's theorem.  Such a fusion of ideas from both sources proved to be particularly crucial in
the recent result \cite{gt-primes} that the primes contained arbitrarily long progressions; this argument was almost entirely finitary in nature, yet at the same time it relied heavily on ideas from the infinitary world of ergodic theory (see \cite{kra-survey}, \cite{host} for further discussion of this connection).

In this paper we investigate a transference in the other direction, taking results from the finitary world of combinatorics (and in particular
graph theory and hypergraph theory), and identifying them with a corresponding result in the infinitary world, which in this case turns out to
be the world of probability theory\footnote{This is actually not all that surprising, given that \emph{finitary} probability theory has already proven to have a major role to play in graph theory.} (or measure theory).  In particular, we present a \emph{correspondence principle}, analogous to
the Furstenberg correspondence principle, that shows how any sequence of increasingly large graphs or hypergraphs has a ``weak limit'', which we view
as an infinitely large \emph{random} graph or hypergraph\footnote{This is related to, but slightly different from, a different concept of graph limit developed by Lov\'asz and Szegedy in \cite{lovasz-szegedy-2}, in which the limiting object becomes a ``continuous weighted graph'', or more precisely a symmetric measurable function from $[0,1] \times [0,1]$ to $[0,1]$.  Such a concrete limiting object is particularly useful for computations such as counting the number of induced subgraphs of a certain shape; it also can be used to establish results such as the triangle removal lemma (Szegedy, personal communication).}.  This principle is slightly more complicated than the Furstenberg correspondence principle, but does not use the full power of deep results such as the Szemer\'edi regularity lemma or its extension to hypergraphs; indeed we do not explicitly state or use such a regularity lemma in this work here, although ideas from that lemma are certainly involved in several components of the argument.  

The main advantage of passing from a deterministic finite graph to a random infinite graph is that one now obtains a number of \emph{factors} ($\sigma$-algebras) in the probability space which enjoy some very useful invariance and relative independence properties. One can think of the presence of these factors as being analogous to the partitions obtained by the Szemer\'edi regularity lemma that make a graph $\eps$-regular, but with the distinction that the partition is now infinite and the $\eps$ parameter set to zero (so one now has perfect regularity).  This sending of the epsilon parameters to zero turns out to be extremely useful in cleaning up proofs of certain statements which previously could only be proven via a regularity lemma.  In particular, we will give an infinitary proof here of the \emph{triangle removal lemma} of Ruzsa and Szemer\'edi \cite{rsz}, as well as the substantially more difficult \emph{hypergraph removal lemma} of Nagle, R\"odl, Schacht, and Skokan \cite{nrs}, \cite{rs}, \cite{rodl}, \cite{rodl2} and Gowers \cite{gowers-hyper} (as well as a later refinement in \cite{tao-hyper}).  As this lemma is already strong enough to deduce Szemer\'edi's theorem on arithmetic progressions (as well as a multidimensional generalisation due to Furstenberg and Katznelson \cite{fk}), we have thus presented yet another proof of Szemer\'edi's theorem here.  These lemmas have some further applications; for instance, they were used in \cite{tao-multiprime} to show that the Gaussian primes contained arbitrarily shaped constellations.  In Appendix \ref{conrecur} we discuss the connections (or lack thereof) between these infinitary removal lemmas, and the recurrence theorems of Furstenberg and later authors.

The setting of this paper was deliberately placed at a midpoint between graph theory and ergodic theory, and the author hopes that it
illuminates the analogies and interconnections between these two subjects.

The author thanks Bal\'asz Szegedy for many useful discussions, Timothy Gowers for suggesting the original topic of investigation, Vitaly Bergelson for encouragement, and Olivier Gerard for corrections.  The author is especially indebted to the anonymous referees for many corrections and suggestions.
The author is supported by a grant from the Packard Foundation.

\section{Motivation: the Furstenberg correspondence principle}\label{motivation}

To motivate the correspondence principle for graphs and hypergraphs, we first review the Furstenberg correspondence principle which
connects results such as Szemer\'edi's theorem with recurrence results in ergodic theory.  Let us recall Szemer\'edi's theorem in a quantitative (finitary) form:

\begin{theorem}[Szemer\'edi's theorem, quantitative version]\label{sz-quant}\cite{szemeredi}  Let $0 < \delta \leq 1$ and $k \geq 1$.  Let $A$ be a subset of a cyclic group $\Z_N := \Z/N\Z$ whose cardinality $|A|$ is at least $\delta N$.  Then there exist at least $c(k,\delta) N^2$ pairs 
$(x,r) \in \Z_N \times \Z_N$ such that $x, x+r, \ldots, x+(k-1)r\in A$, where $c(k,\delta) > 0$ is a positive quantity depending only on $k$ and $\delta$.
\end{theorem}

This result is easily seen to imply to Szemer\'edi's theorem in its traditional (infinitary) form, which asserts that every set of integers of positive upper density contains arbitrarily long progressions.  The converse implication also follows from an argument of Varnavides \cite{var}.  This particular formulation of Szemer\'edi's theorem played an important role in the recent result \cite{gt-primes} that the primes contained arbitrarily long arithmetic progressions.

In 1977, Furstenberg obtained a new proof of Szemer\'edi's theorem by deducing it from the following result in ergodic theory.

\begin{theorem}[Furstenberg recurrence theorem]\label{furst-thm}\cite{furst}, \cite{furstenberg}  Let $(\Omega, \B_{\max}, \P)$ be a probability space (see Appendix \ref{probtools} for probabilistic notation).  Let $T:\Omega \to \Omega$ be a bi-measurable map which is probability preserving, thus $\P( T^n A ) = \P(A)$ for all events $A \in \B_{\max}$ and $n \in \Z$.  Then for all $k \geq 1$ and all events $A \in \B_{\max}$ with $\P(A) > 0$, we have
$$ \liminf_{N \to \infty} \frac{1}{N} \sum_{n=1}^N \P( A \wedge T^n A \wedge \ldots \wedge T^{(k-1)n} A ) > 0.$$
\end{theorem}

The deduction of Theorem \ref{sz-quant} from Theorem \ref{furst-thm} proceeds by the \emph{Furstenberg correspondence principle} \cite{furst}, \cite{furstenberg}, \cite{furst-book}.  Let us give a slightly non-standard exposition of this principle (in particular 
drawing heavily on the language of probability theory), in order to motivate an 
analogous principle for graphs and hypergraphs in later sections.
We shall interpret this correspondence principle as an assertion that any sequence $(A^{(m)}, \Z_{N^{(m)}})$ of sets $A^{(m)}$ in a cyclic group
$\Z_{N^{(m)}}$ can have an asymptotic limit as $m \to \infty$, which will end up being a probability space endowed with a probability-preserving shift $T$.  To state this more precisely we shall need some notation.  First, we describe a certain universal space in which it will be convenient to
take limits.

\begin{definition}[Furstenberg universal space]
Let $\Omega := 2^\Z := \{ B: B \subset \Z\}$ denote the set of all subsets $B \subset \Z$ of the integers $\Z$; one can also view this space as the infinite cube $\{0,1\}^\Z$ if desired.  We give this space the product $\sigma$-algebra $\B_{\max}$, generated by the events\footnote{A more topological of thinking about this proceeds by endowing $\Omega$ with the product topology, so that it becomes a compact Hausdorff totally disconnected space, and then letting $\B_{\max}$ be the Borel $\sigma$-algebra, generated by the open sets.  The regular algebra $\B_\reg$ then consists those events which are simultaneously open and closed, or equivalently those events whose indicator function is continuous.} $A_n := \{ B \in \Omega: n \in B \}$ for $n \in \Z$.  Indeed, one can think of $(\Omega,\B_{\max})$ as being the universal event space generated by the countable sequence of events $A_n$.  The space $\Omega$ enjoys an obvious shift action $T: \Omega \to \Omega$, defined by $TB := B+1 :=\{ n+1: n \in B\}$ for all $B \in \Omega$.  This then induces a shift $T: \B_\max \to \B_\max$ in the obvious manner, thus for instance $T^n A_0 = A_n$.  We define the \emph{regular algebra} $\B_\reg$ of $\B_{\max}$ to be the algebra generated by the $A_n$, thus the events in $\B_\reg$ (which we refer to as \emph{regular events}) are those events which are generated by at most finitely many of the $A_n$ (i.e. those events that only require knowing the truth value of $n \in B$ for finitely many values of $n$).  
\end{definition}

Now we embed finite objects $(A^{(m)}, \Z_{N^{(m)}})$ described earlier in this universal space.

\begin{definition}[Furstenberg universal embedding]\label{furst-embed}
Let $m \geq 1$, let $\Z_{N^{(m)}}$ be a cyclic group with $N^{(m)} \geq m$, and let $A^{(m)}$ be a subset of $\Z_{N^{(m)}}$.
We define the probability space $(\Omega^{(m)}, \B_{\max}^{(m)}, \P^{(m)})$ as the space corresponding to sampling\footnote{The introduction of the dilation parameter $\lambda^{(m)}$ is essentially the averaging trick of Varnavides \cite{var}. The exact construction of this space is not important so long as one has the independent random variables $x^{(m)}$ and $\lambda^{(m)}$, but for sake of concreteness one can set $\Omega^{(m)} := \Z_{N^{(m)}} \times [L^{(m)}]$, $\B_{\max}^{(m)} := 2^{\Omega^{(m)}}$ to be the power set of $\Omega^{(m)}$, and $\P^{(m)}$ to be the uniform distribution on $\Omega^{(m)}$.} $x^{(m)}$ and $\lambda^{(m)}$ uniformly and independently at random from $\Z_{N^{(m)}}$ and $[L^{(m)}]$, where $[N] := \{1,\ldots,N\}$ denotes the integers from $1$ to $N$ and 
$L^{(m)} \geq 1$ is the integer part of $N^{(m)}/m$.  We then map every pair $(x^{(m)}, \lambda^{(m)})$ of $\Omega^{(m)}$
to a point $B^{(m)} \in \Omega$ (i.e. a subset of the integers $\Z$) by the formula
$$ B^{(m)} := \{ n \in \Z: x^{(m)} + n \lambda^{(m)} \in A^{(m)} \};$$
one can think of this as a random lifting of the set $A^{(m)} \subset \Z_{N^{(m)}}$ up to the integers $\Z$.  This mapping from $\Omega^{(m)}$
to $\Omega$ is clearly measurable, since the inverse images of the generating events $A_n$ in $(\Omega, \B_{\max})$ are simply the events
that $x^{(m)} + n \lambda^{(m)} \in A^{(m)}$, which are certainly measurable in $\B^{(m)}_\max$.  This allows us to extend the probability measure
$\P^{(m)}$ from $(\B_{\max}^{(m)}, \Omega^{(m)})$ to the product space $(\B_{\max} \times \B_{\max}^{(m)}, \Omega \times \Omega^{(m)})$ 
in a canonical manner\footnote{More precisely, we graph the measurable mapping from $\Omega^{(m)}$ to $\Omega$ as a measurable mapping from $\Omega^{(m)}$ to $\Omega \times \Omega^{(m)}$, which contravariantly induces a $\sigma$-algebra homomorphism from the product $\sigma$-algebra $\B_\max \times \B_\max^{(m)}$ to $\B_\max^{(m)}$.  Pulling back the probability measure $\P^{(m)}$ under this homomorphism yields the extension.  A similar construction applies to the graph and hypergraph embeddings in later sections.}, identifying the events $A_n$ with the events $x^{(m)} + n \lambda^{(m)} \in A^{(m)}$.  We shall abuse notation and
refer to the extended measure also as $\P^{(m)}$.
\end{definition}

In more informal terms, the Furstenberg embedding has created, for each $m$, a random set $B^{(m)} \subset \Z$ which
will capture all the important information about the original set $A^{(m)}$ and $\Z_{N^{(m)}}$.  For instance, the density of $A^{(m)}$ is
nothing more than the probability that $0$ lies in $B^{(m)}$, or equivalently the probability of the event $A_0$.  One can view $x^{(m)}$ and
$\lambda^{(m)}$ as the ``hidden variables'' which generate this random set $B^{(m)}$.  However, in order to invoke the correspondence principle we will need to ``forget'' that the random set $B^{(m)}$ actually came from these variables; indeed, we are going to restrict $\P^{(m)}$ to the
common factor $(\Omega, \B_{\max})$ in order to take limits as $m \to \infty$.  More precisely, we have

\begin{proposition}[Furstenberg correspondence principle]  For every $m \geq 1$, let $\Z_{N^{(m)}}$ be a cyclic group with $N^{(m)} \geq m$, and let $A^{(m)}$ be a subset of $\Z_{N^{(m)}}$, and let $\P^{(m)}$ be as in Definition \ref{furst-embed}.  Then there exists a subsequence $0 < m_1 < m_2 < \ldots$ of $m$, and a probability measure $\P^{(\infty)}$ on the Furstenberg universal space $(\Omega,\B_{\max})$, such that we have the weak convergence property
\begin{equation}\label{weak-conv}
\lim_{i \to \infty} \P^{(m_i)}( E ) = \P^{(\infty)}( E ) \hbox{ for all } E \in \B_{\reg}.
\end{equation}
Furthermore, we have the shift invariance property
\begin{equation}\label{pshift}
\P^{(\infty)}( T^n E ) = \P^{(\infty)}(E) \hbox{ for all } E \in \B_{\max}, n \in \Z.
\end{equation}
\end{proposition}

\begin{proof}  The algebra $\B_{\reg}$ is countable.  Thus the existence of the weak limit $\P^{(m_i)}$ follows from Lemma \ref{arzela}.
Now observe that the random sets $B^{(m)}$ and $T^n B^{(m)} = B^{(m)}+n$ have the same probability distribution (because $x_0$ and $x_0-n\lambda$
have the same distribution for any fixed $\lambda$).  Thus we observe that $\P^{(m)}$ is shift-invariant:
$$
\P^{(m)}( T^n E ) = \P^{(m)}(E) \hbox{ for all } E \in \B_\max, n \in \Z.
$$
Applying \eqref{weak-conv} we obtain \eqref{pshift} for all \emph{regular} events $E$.  But since $\P^{(\infty)}$ is countably additive, we see that the space of events $E$ for which \eqref{pshift} holds for every $T$ is a $\sigma$-algebra which contains $\B_\reg$, and thus contains $\B_\max$ as claimed.
\end{proof}

Now we can deduce Theorem \ref{sz-quant} from Theorem \ref{furst-thm}.

\begin{proof}[of Theorem \ref{sz-quant} assuming Theorem \ref{furst-thm}]
Suppose that Theorem \ref{sz-quant} fails.  Then we can find $k \geq 1$ and $0 < \delta \leq 1$, a sequence $N^{(m)}$ of positive integers, and a sequence of sets $A^{(m)} \subset \Z_{N^{(m)}}$ of density $|A^{(m)}|/N^{(m)} \geq \delta$ such that
$$ \lim_{m \to \infty} \frac{1}{(N^{(m)})^2} | \{ (x,r) \in \Z_{N^{(m)}}: x, x+r, \ldots, x+(k-1)r \in A^{(m)} \} | = 0.$$
By passing to a subsequence of $m$ if desired we can make this convergence arbitrarily fast; for instance, we can ensure that
\begin{equation}\label{nmform}
\frac{1}{(N^{(m)})^2} | \{ (x,r) \in \Z_{N^{(m)}}: x, x+r, \ldots, x+(k-1)r \in A^{(m)} \} | \leq \delta 100^{-m}.
\end{equation}
Observe that the left-hand side is at least $\frac{1}{(N^{(m)})^2} |A^{(m)}| \geq \delta / N^{(m)}$, so we conclude that
$$ N^{(m)} \geq 100^m.$$
In particular $N^{(m)} \geq m$, so we can invoke the Furstenberg correspondence principle and obtain a shift-invariant system
$(\Omega, \B_{\max}, \P^{(\infty)})$ on the Furstenberg universal space $(\Omega, \B_{\max})$ with the stated properties.

Now let us compute some probabilities in this system, starting with the probability of $A_0$.  From definition of $\P^{(m)}$ we have
$$
\P^{(m)}( A_0 ) = \P^{(m)}( x_0 \in A^{(m)} ) 
= |A^{(m)}| / N^{(m)}
\geq \delta
$$
so by \eqref{weak-conv} we have
$$ \P(A_0) \geq \delta.$$
In particular $A_0$ has strictly positive probability.  Next, let $1 \leq n \leq m$ and consider the expression
\begin{align*}
\P^{(m)}( A_0 \wedge T^n A_0 \wedge \ldots \wedge T^{(k-1)n} A_0) &= \P^{(m)}(A_0 \wedge A_n \wedge \ldots \wedge A_{(k-1)n} ) \\
&= \P^{(m)}( x^{(m)} + j n \lambda^{(m)} \in A^{(m)}\ \forall 0 \leq j < k) \\
&= \frac{1}{N^{(m)} L^{(m)}} | \{ (x^{(m)}, \lambda^{(m)}) \in \Z_{N^{(m)}} \times [L^{(m)}]: \\
&\quad x^{(m)}+jn\lambda^{(m)} \in A^{(m)}\ \forall 0 \leq j < k \}|.
\end{align*}
Now observe from definition of $L^{(m)}$, the progressions $x^{(m)}, x^{(m)}+n\lambda^{(m)}, \ldots, x^{(m)} + (k-1) n \lambda^{(m)}$ are all distinct as $x^{(m)}$ and $\lambda^{(m)}$
vary.  Applying \eqref{nmform}, we see that
\begin{align*}
\P^{(m)}( A_0 \wedge T^n A_0 \wedge \ldots \wedge T^{(k-1)n} A_0) &= \P^{(m)}(A_0 \wedge A_n \wedge \ldots \wedge A_{(k-1)n} ) \\
&\leq \frac{1}{N^{(m)} L^{(m)}} (N^{(m)})^2 \delta 100^{-m} \\
&\leq 2 m 100^{-m}
\end{align*}
(say) for all $1 \leq n \leq m$.  In particular we have
$$
\lim_{m \to \infty} \P^{(m)}( A_0 \wedge T^n A_0 \wedge \ldots \wedge T^{(k-1)n} A_0) = 0
$$
for each fixed $n \geq 1$, and hence by \eqref{weak-conv}
$$ \P^{(\infty)}( A_0 \wedge T^n A_0 \wedge \ldots \wedge T^{(k-1)n} A_0 ) = 0$$
for all $n \geq 1$.  But this contradicts Theorem \ref{furst-thm}.  This completes
the deduction of Theorem \ref{sz-quant} from Theorem \ref{furst-thm}. 
\end{proof}

\begin{remark} Note that as this proof proceeded by contradiction, it does not obviously give any sort of \emph{quantitative} lower bound
for the quantity $c(k,\delta)$ appearing in Theorem \ref{sz-quant}.  It is actually possible (with nontrivial effort) to extract such a bound by taking the \emph{proof} of Theorem \ref{furst-thm} and making everything finitary; see \cite{tao:ergodic}.  However the bounds obtained in this manner are extremely poor.  The same remarks apply to the infinitary proofs of the triangle removal lemma and hypergraph removal lemma that we give below.  As a related remark, observe that the above argument, while infinitary, did not require the axiom of choice, as one can eliminate the apparent use of choice at the beginning of the argument by well-ordering the objects $A$, $\Z_N$, $\delta$ appearing in Theorem \ref{sz-quant} in some standard manner.  (The use of Lemma \ref{arzela} also does not require choice; see Remark \ref{construct}.  The original proof of the Furstenberg recurrence theorem in \cite{furst} is also choice-free, though the later proof in \cite{furstenberg} is not, as it uses Zorn's lemma.) Indeed we will not actually need the axiom of choice in this entire paper, though we shall assume it in order to simplify the exposition slightly.
\end{remark}

\begin{remark} One can also reverse the above argument, and use Theorem \ref{sz-quant} to deduce Theorem \ref{furst-thm}, basically by applying Theorem \ref{sz-quant} to various truncated versions of the random set $B := \{ n \in \Z: T^n x \in A \}$, where $x$ is sampled from the sample space $\Omega$ using the probability measure $\P$.  We omit the standard details.
\end{remark}

\section{The graph correspondence principle}\label{graphsec}

We now develop an analogue of the Furstenberg correspondence principle for graphs; namely, we start with a sequence of (undirected) graphs
$G^{(m)} = (V^{(m)}, E^{(m)})$ for each $m \geq 1$, and wish to extract (after passing to a subsequence of $m$'s) some sort of infinitary weak limit.
This type of problem was already addressed in \cite{lovasz-szegedy-2}, with the main tool being a certain weak form of
the Szemer\'edi regularity lemma.  Our approach is somewhat similar (though not identical), and the regularity lemma will appear only after the infinite limit is extracted, in Lemma \ref{irl} below.

As before, we need a universal space in which to take limits.  Just as the Furstenberg universal space consisted of infinite \emph{sets} of integers, the graph universal space will consist of infinite \emph{graphs} on the natural numbers.  The shift $T$ (which represents a $\Z$-action) is now replaced\footnote{We are indebted to Bal\'asz Szegedy for pointing out the analogy between the $\Z$-action of a dynamical system and the $S_\infty$-action on an infinite graph.} by the action of the permutation group $S_\infty$, defined
as the group of all permutations $\sigma: \Z \to \Z$ of the integers.

\begin{definition}[Graph universal space]\label{graph-universal}
Let $\N := \{1, 2, \ldots\}$ denote the natural numbers, and let
$\Omega := 2^{\binom{\N}{2}} = \{ (\N, E_\infty): E_\infty \subset \binom{\N}{2} \}$ denote the space of all (infinite) graphs $(\N, E_\infty)$ on the natural numbers, thus the edge set $E_\infty$ is an arbitrary collection of unordered pairs of distinct integers.  On this space $\Omega$, we introduce the events $A_{i,j} = A_{j,i}$ for any unordered pair of distinct natural numbers $\{i,j\} \in \binom{\N}{2}$ by $A_{i,j} := \{ (\N,E_\infty) \in \Omega: (i,j) \in E_\infty \}$, and let $\B_{\max}$ be the $\sigma$-algebra generated by the countable sequence of events $A_{i,j}$.  (We adopt the convention that $A_{i,i} = \emptyset$ for all $i \in \N$, thus our graphs have no loops.)  
We also introduce the regular algebra $\B_{\reg}$ generated by the $A_{i,j}$, thus these are the events that depend only on finitely many of 
the $A_{i,j}$.
For any permutation $\sigma: \N \to \N$ of the natural numbers, we define the associated action on $\B_{\max}$ by mapping $\sigma: A_{i,j} \mapsto A_{\sigma(i),\sigma(j)}$ and extending this to a $\sigma$-algebra isomorphism in the unique manner.  More explicitly, $\sigma$ will map each graph $(N, E_\infty)$ to the graph $(N, \sigma E_\infty)$, where $\sigma E_\infty := \{ \{\sigma(i), \sigma(j)\}: \{i,j\} \in E_\infty \}$.
For any (possibly infinite) subset $I$ of $\N$, we define $\B_I$ to be the factor of $\B_{\max}$ generated by
the events $A_{i,j}$ for $i,j \in I$; informally speaking, $\B_I$ represents the knowledge obtained by measuring the restriction of $E_\infty$ to $I$.
Observe the trivial monotonicity $\B_I \subseteq \B_J$ whenever $I \subseteq J$.
\end{definition}

The space $(\Omega, \B_{\max})$ is thus the universal event space associated to the events $A_{i,j}$, and is the natural event space for studying infinite random graphs. (For instance, the infinite Erd\"os-Renyi random graph $G(\infty,p)$ for fixed $0 \leq p \leq 1$, where the vertex set is $\Z$ (say) and any two integers are connected by an edge with an independent probability of $p$, would correspond to the scenario in which all the events $A_{i,j}$ are independent with probability $p$ each.)  The permutation group $S_\infty$ defined earlier acts on the event space $(\Omega, \B_{\max})$ in the obvious manner.   Thus for instance $\sigma(\B_I)=\B_{\sigma(I)}$ for all $\sigma \in S_\infty$ and $I \subseteq \N$.

Next, we need a way to embed every finite graph into the universal space.

\begin{definition}[Graph universal embedding]\label{graph-embed}  Let $m \geq 1$, and let $G^{(m)} = (V^{(m)}, E^{(m)})$ be a finite graph.
Let $(\Omega^{(m)}, \B_{\max}^{(m)}, \P^{(m)})$ be the probability space corresponding to the sampling of a countable sequence\footnote{This sequence contains the ``hidden variables'' that will play the role of the parameters $x^{(m)}$ and $\lambda^{(m)}$ in the preceding section.  Again, the exact construction of this Wiener-type probability space is not important. The most canonical way to proceed is to let $\Omega^{(m)}$ be the countable product
$(V^{(m)})^\N$ with the product $\sigma$-algebra $\B_{\max}^{(m)}$ and the product uniform probability measure $\P^{(m)}$.  A more concrete way
would be to identify $V^{(m)}$ with $[n^{(m)}] = \{1,\ldots,n^{(m)}\}$ by appropriate labeling, 
set $\Omega^{(m)}$ to be the unit interval $[0,1) := \{ x: 0 \leq x < 1 \}$, let $\B_{\max}^{(m)}$ be the Borel $\sigma$-algebra,
$\P^{(m)}$ be Lebesgue measure, and let $x_j$ be the $j^{th}$ digit in the base-$n^{(m)}$ expansion of $x$ (rounding down when a terminating decimal occurs).} of i.i.d. random variables $x^{(m)}_1, x^{(m)}_2, \ldots \in V^{(m)}$ sampled independently and uniformly at random\footnote{Of course for any fixed $m$ there will be infinitely many repetitions among these $x^{(m)}_i$ since $V^{(m)}$ is finite, but in practice we are interested in taking limits in which $|V^{(m)}| \to \infty$, and so these collisions will become asymptotically negligible.}.  To every sequence
$(x^{(m)}_1,x^{(m)}_2, \ldots ) \in \Omega^{(m)}$ we associate an infinite graph $G^{(m)}_\infty = (\N, E^{(m)}_\infty) \in \Omega$ by setting
$$ E^{(m)}_\infty := \{ \{i,j\} \in \binom{\N}{2}: \{x^{(m)}_i, x^{(m)}_j\} \in E^{(m)} \};$$
one can think of this as a random lifting of the graph $G^{(m)}$ on $V^{(m)}$ up to an infinite graph $G^{(m)}_\infty$ on the natural numbers $\N$.
This mapping from $\Omega^{(m)}$
to $\Omega$ is clearly measurable, since the inverse images of the generating events $A_{i,j}$ in $(\Omega, \B_{\max})$ are simply the events
that $\{x^{(m)}_i, x^{(m)}_j\}$ lie in $G^{(m)}$, which are certainly measurable in $\B^{(m)}_\max$.  
This allows us to extend the probability measure
$\P^{(m)}$ from $(\B_{\max}^{(m)}, \Omega^{(m)})$ to the product space $(\B_{\max} \times \B_{\max}^{(m)}, \Omega \times \Omega^{(m)})$ 
in a canonical manner, identifying the events $A_{i,j}$ with the events $\{ x^{(m)}_i, x^{(m)}_j \} \in E^{(m)}$.  We shall abuse notation and
refer to the extended measure also as $\P^{(m)}$.
\end{definition}

\begin{remarks} Now that the space $(\Omega^{(m)}, \B_{\max}^{(m)})$ is infinite, not every event involving the $x^{(m)}_i$ is measurable, however any event which involves only finitely many of the $x^{(m)}_i$ is certainly measurable (and in particular has a well-defined probability).
One can view $E^{(m)}_\infty$ as the infinite random graph formed by statistically sampling of the original finite (and deterministic) graph $E^{(m)}$.  This is a convenient way to convert arbitrary graphs, on arbitrary vertex sets, to a fixed universal (random) graph on a fixed universal vertex set, in this case the natural numbers $\N$.  The random graph $E^{(m)}_\infty$ turns out to capture all the relevant features we require of the original graph; for instance, the probability that $E^{(m)}_\infty$ lies in the event $A_{1,2}$ is essentially\footnote{We say ``essentially'' because there is a slight error term coming from the event that $x^{(m)}_1=x^{(m)}_2$.  However this error will become negligible in limits for which $|V^{(m)}| \to \infty$.} the edge density of $E^{(m)}$, while the probability that $E^{(m)}_\infty$ lies in $A_{1,2} \wedge A_{2,3} \wedge A_{3,1}$ is essentially the triangle density of $E^{(m)}$, and so forth.  On the other hand, it suppresses irrelevant features such as what the labels of the original vertex set $V^{(m)}$ were; in particular, applying a graph isomorphism to $E^{(m)}$ does not affect the probability distribution of $E^{(m)}_\infty$ at all.  More generally, we observe the permutation invariance
\begin{equation}\label{perm}
\P^{(m)}( \sigma E ) = \P^{(m)}(E) \hbox{ for all } E \in \B_{\max}, \sigma \in S_\infty
\end{equation}
which can be verified by first checking on regular events $E$ (i.e. finite boolean combinations of the $A_{i,j}$) and then extending as
in the proof of the Furstenberg correspondence principle.
\end{remarks}

Once again, we can view the random graph $G^{(m)}_\infty$ as being generated\footnote{This is of course the perspective taken in property testing.  It is not surprising that the Szemer\'edi regularity lemma plays a crucial role in that theory also; see \cite{alon-shapira}.  Indeed, this argument suggests that an infinitary approach to property testing theory is possible.}
 by ``hidden variables'' $x^{(m)}_1, x^{(m)}_2, \ldots$.  As before, we wish to ``forget'' these hidden variables and pass to a limit.  This
can be achieved as follows.

\begin{proposition}[Graph correspondence principle]\label{gcp}  For every $m \geq 1$, let $G^{(m)} = (V^{(m)}, E^{(m)})$ be a finite undirected graph, and let $\P^{(m)}$ be as in Definition \ref{graph-embed}.  Then there exists a subsequence $0 < m_1 < m_2 < \ldots$ of $m$, and a probability measure $\P^{(\infty)}$ on the graph universal space $(\Omega,\B_{\max})$, such that we have the weak convergence property \eqref{weak-conv}
and the permutation invariance property
\begin{equation}\label{pshift-2}
\P^{(\infty)}( \sigma E ) = \P^{(\infty)}(E) \hbox{ for all } E \in \B_{\max}, \sigma \in S_\infty.
\end{equation}
\end{proposition}

\begin{proof} The algebra $\B_{\reg}$ is countable.  Thus the existence of the weak limit $\P^{(m_i)}$ follows from Lemma \ref{arzela}.
From \eqref{perm} we can deduce \eqref{pshift-2} by arguing exactly as in the Furstenberg correspondence principle.
\end{proof}

So far, the permutation group $S_\infty$ has played the same role for graphs as the integer group $\Z$ played for sets of integers.  However, the permutation group is significantly more ``mixing'', which allows us to immediately ``regularise'' the system obtained in Proposition \ref{gcp}:

\begin{lemma}[Infinitary regularity lemma]\label{irl}  
Let $\P$ be a probability measure on the graph universal space $(\Omega, \B_\max)$ which is permutation-invariant in the sense of \eqref{pshift-2}.  Then for any 
$I, I_1,\ldots, I_l \subset \N$ with $I \cap I_1 \cap \ldots \cap I_l$ infinite, the factors $\B_I$ and $\bigvee_{i=1}^l \B_{I_i}$ are relatively independent conditioning on $\bigvee_{i=1}^l \B_{I \cap I_i}$, with respect to this probability measure $\P$.  
(See Appendix \ref{probtools} for a definition of relative independence.)
\end{lemma}

This result is the infinitary analogue of the Szemer\'edi regularity lemma, and will play a crucial role in establishing the proof of the triangle removal lemma (and later, the hypergraph removal lemma) in subsequent sections.  

\begin{proof} Fix $I, I_1, \ldots, I_l$.  We may assume $l \geq 1$ 
since the claim is trivial when $l=0$.
To show that $\B_I$ and $\bigvee_{i=1}^l \B_{I_i}$ are relatively independent conditioning on
$\bigvee_{i=1}^l \B_{I \cap I_i}$ with respect to the probability measure $\P^{(\infty)}$, it suffices by Lemma \ref{verify-independence} to show that
$$ \| \P( E_I | \bigvee_{i=1}^l \B_{I_i} ) \|_{L^2} = \| \P( E_I | \bigvee_{i=1}^l \B_{I \cap I_i} ) \|_{L^2}$$
for all $E_I \in \B_I$.  By Lemma \ref{dense} and limiting arguments we may assume without loss of generality that $E_I$ is regular.  In particular we have $E_I \in \B_{I'}$ for some finite subset $I'$ of $I$.  
By Corollary \ref{product-monotone} and a limiting argument we may assume that the set $I$ has an 
infinite complement.  By another such limiting argument we can also assume that $I_i \backslash I$ is finite for all $i$.

Let $A$ be the infinite set $A := (I \cap I_1 \cap \ldots \cap I_l) \backslash I'$, and let $B$ be the finite set $B := \bigcup_{i=1}^l I_i \backslash I$.  Then we can find a permutation $\sigma$ be a permutation which maps $A$ to $A \cup B$ bijectively, but is constant on $I \backslash A$, and in particular fixes $I'$.  Thus $\sigma$ also fixes $E_I$, and maps $I \cap I_i$ to $(I \cap I_i) \cup B$.  Thus
$$ \| \P( E_I | \bigvee_{i=1}^l \B_{(I \cap I_i) \cup B} ) \|_{L^2} = \| \P( E_I | \bigvee_{i=1}^l \B_{I \cap I_i} ) \|_{L^2}.$$
But as $\B_{I_i}$ lies between $\B_{I \cap I_i}$ and $\B_{(I \cap I_i) \cup B}$, the claim now follows from 
Lemma \ref{pythag}.  
\end{proof}

\begin{remark} The above proof of the regularity lemma is short but perhaps a bit opaque.  Let us informally discuss a special case of this lemma, namely that the events $A_{1,3}$ and $A_{2,3}$ are relatively independent conditioning on $\B_{\{3,4,5,\ldots\}}$; this is a special case of the situation where $I = \N \backslash \{2\}$, $I_1 = \N \backslash \{1\}$, and $l=1$.  
Passing back to the finite graph setting (by working with the probability measures $\P^{(m)}$ from Proposition \ref{gcp}), this claim may seem puzzling at first, because the
events $\{x^{(m)}_1, x^{(m)}_3\} \in E^{(m)}$ and  $\{x^{(m)}_2, x^{(m)}_3\} \in E^{(m)}$ can certainly be correlated; indeed, whenever $x^{(m)}_3$
has high degree, then both events occur with high probability, and when it has low degree, both events occur with low probability.  However, if one can somehow learn the degree of $x^{(m)}_3$, then these two events become relatively independent \emph{conditioning on the degree of $x^{(m)}_3$}.
And now the purpose of the factor $\B_{\{3,4,5,\ldots\}}$ becomes clear; by ``polling'' many additional 
vertices $x^{(m)}_4, x^{(m)}_5, \ldots, x^{(m)}_N$ and measuring the connectivity of $x^{(m)}_3$ with all of these additional vertices, we can obtain
a statistical prediction for the degree of $x^{(m)}_3$, whose accuracy and confidence level become almost surely perfect in the asymptotic limit $N \to \infty$.
More generally, it turns out that by polling the interconnectivity of vertices in the infinite set $I \cap I_i$ for $i=1,\ldots,l$ one can obtain an almost surely perfectly accurate prediction of all the ``common information'' held between an event in $\B_I$ and an event in $\bigvee_{i = 1}^l \B_{I_i}$.  Let us illustrate this with one further example, namely the relative independence of $A_{1,2}$ and $A_{2,3} \wedge A_{1,3}$
conditioning on $\B_{\{1,4,5,\ldots\}} \vee \B_{\{2,4,5,\ldots\}}$; this corresponds to the case $I = \N \backslash \{3\}$, $l=2$, and
$I_i = \N \backslash \{i\}$ for $i=1,2$.  We are asking for the events $\{x^{(m)}_1, x^{(m)}_2 \} \in E^{(m)}$ and
$\{x^{(m)}_1, x^{(m)}_3 \}, \{x^{(m)}_2, x^{(m)}_3 \} \in E^{(m)}$ to become relatively independent once we sample all the connectivity information
between $x^{(m)}_2$ and $x^{(m)}_4, x^{(m)}_5, \ldots$, and between $x^{(m)}_3$ and $x^{(m)}_4, x^{(m)}_5, \ldots$.  To see how this will work,
observe that while the two events in question will not be unconditionally independent in general, they will become conditionally independent once
the number of paths of length two connecting $x^{(m)}_1$ and $x^{(m)}_2$ are known, since upon freezing $x^{(m)}_1$ and $x^{(m)}_2$ this determines
the probability that the independent variable $x^{(m)}_3$ will satisfy the latter 
event $\{x^{(m)}_1, x^{(m)}_3 \}, \{x^{(m)}_2, x^{(m)}_3 \} \in E^{(m)}$; since $x^{(m)}_3$ does not affect the former event
 $\{x^{(m)}_1, x^{(m)}_2 \} \in E^{(m)}$, we obtain relative independence.  But the number of paths of length two can be determined statistically,
 by counting the proportion of $j \in \{4,5,\ldots\}$ for which $\{ x^{(m)}_1, x^{(m)}_j \}$ and $\{ x^{(m)}_2, x^{(m)}_j \}$
both lie in $E^{(m)}$.  This lies in the factor $\B_{\{1,4,5,\ldots\}} \vee \B_{\{2,4,5,\ldots\}}$ and is the reason for the conditional
 independence\footnote{There is another way of viewing this, namely that each vertex $x^{(m)}_j$ induces a partition of the $x^{(m)}_1$ and $x^{(m)}_2$ vertex sets, by dividing them into those vertices which are connected to $x^{(m)}_j$ in $E^{(m)}$ and those that are not.  Letting $j$ vary in $\{4,5,\ldots,N\}$ one obtains a partition of these vertex classes which behaves increasingly like the partitions created by the Szemer\'edi regularity lemma as $N \to \infty$, in the sense that the graph between the $x^{(m)}_1$ and $x^{(m)}_2$ becomes increasingly ``$\eps$-regular'' relative to this partition; the $\eps$-regularity is closely related to the relative independence properties discussed here.  We will however not pursue this approach as it becomes somewhat complicated when we move to the hypergraph setting, whereas the techniques we present here carries over to hypergraphs with virtually no changes.}.
\end{remark}

\begin{remark} Similar correspondence principles exist for bipartite graphs, directed graphs, and multicolored graphs (where the color set is kept independent of $m$), and so forth; for instance, a generalisation to tripartite graphs is sketched out in Appendix \ref{conrecur}. We will not pursue the other generalisations here as they are rather minor, though we will consider a hypergraph extension of this principle in Section \ref{hyperspace}.
\end{remark}

\section{An infinitary proof of the triangle removal lemma}\label{triangle-sec}

Let us now apply the above correspondence principle to obtain the following triangle-removal lemma of Ruzsa and Szemer\'edi:

\begin{lemma}[Triangle removal lemma]\label{trl}\cite{rsz} Let $G = (V,E)$ be an undirected graph with $|V|=n$ vertices.  Suppose that $G$ contains fewer than $\delta n^3$ triangles for some $0 < \delta \leq 1$, or more precisely
$$ | \{ (x_1,x_2,x_3) \in V^3: \{x_1,x_2\}, \{x_2,x_3\}, \{x_3,x_1\} \in E \} | \leq \delta n^3.$$
Then it is possible to delete $o_{\delta \to 0}(n^2)$ edges from $G$ to create a graph $G'$ which is triangle-free.  Here $o_{\delta \to 0}(n^2)$ denotes a quantity, which when divided by $n^2$, goes to zero as $\delta \to 0$, uniformly in $n$.
\end{lemma}

Previous to this paper, the only known proof of lemma proceeded via the \emph{Szemer\'edi regularity lemma} \cite{szemeredi}.  It can be used among other things to imply the $k=3$ case of Szemer\'edi's theorem (Theorem \ref{sz-quant}).  Based on this connection, it is natural to ask whether there is an infinitary analogue of this lemma, similarly to how Theorem \ref{furst-thm} is an infinitary counterpart to Theorem \ref{sz-quant}.  We shall deduce it from the following substantially stronger infinitary statement.

If $J$ is a set, we define an \emph{downset} $\ideal$ in $J$ to be any collection of subsets $e$ of $J$ with the property that whenever $e \in \ideal$ and $e' \subseteq e$, then $e' \in \ideal$ also.  In particular, downsets are automatically closed under intersection.

\begin{theorem}[Hypergraph removal lemma, infinitary version]\label{compact-triangle}  Let $(\Omega, \B_{\max}, \P)$ be a probability space, and let $\B_{\reg} \subseteq \B_{\max}$ be an algebra.  Let $J$ be a finite set, and let $\ideal_\max$ be an downset in $J$. 
For each $e \in \ideal_\max$ let $\B_e$ be a factor of $\B_{\max}$ with
the following properties:
\begin{itemize}
\item (Regularisability) Each of the factors $\B_e$ is generated by countably many events from $\B_\reg$.
\item (Nesting) If $e, e' \in \ideal_\max$ are such that $e \subseteq e'$, then $\B_e$ is a factor of $\B_{e'}$.
\item (Independence) If $e, e_1,\ldots, e_l \in \ideal_\max$, then the factors $\B_e$ and $\bigvee_{i=1}^l \B_{e_i}$ are relatively independent conditioning on $\bigvee_{i=1}^l \B_{e \cap e_i}$.
\end{itemize}
For each $e \in \ideal_\max$, let $E_e$ be an event in $\B_e$ such that
$$ \P( \bigwedge_{e \in \ideal_\max} E_e ) = 0.$$
Then for any $\eps > 0$, there exist events $F_e \in \B_e \cap \B_{\reg}$ for all $e \in \ideal_\max$
such that
$$\P( E_e \backslash F_e ) \leq \eps \hbox{ for all } e \in \ideal_\max$$
and
$$ \bigwedge_{e \in \ideal_\max} F_e = \emptyset.$$
\end{theorem}

We will prove this rather strange-looking proposition in Section \ref{hypersec}.  For the purposes of proving the triangle removal lemma, we will only need this lemma in the special case when $J = \{1,2,3\}$, when $\ideal_\max := \{ e: |e| \leq 2 \}$, and when $A_e = \Omega$ for all $e \neq \{1,2\}, \{2,3\}, \{3,1\}$.  However, the lemma is not that much more difficult to prove in the general case\footnote{This is in stark contrast to the finitary situation, in which the hypergraph removal lemma is significantly more difficult than the triangle removal lemma.  This is ultimately because of the need in the finitary hypergraph setting to constantly play off epsilons of different sizes against each other; see \cite{nrs}, \cite{rs}, \cite{rodl}, \cite{rodl2}, \cite{gowers-hyper}, \cite{tao-hyper} for some examples of this.  However in the infinitary asymptotic limit, most of the epsilons have disappeared or at least been confined to individual lemmas where they do not interact with other epsilons.  This simplifies the proof significantly, albeit at the cost of working in an infinitary setting as opposed to a finitary one.  In the converse direction, note the proliferation of epsilons in \cite{tao:ergodic} when Furstenberg's proof of Szemer\'edi's theorem is transferred from the infinitary setting to the finitary one.}, and it will rather easily yield a hypergraph generalisation of the triangle removal lemma, so we retain the proposition in the general form.
The hypothesis $\P( \bigwedge_{e \in \ideal_\max} E_e ) = 0$ is the analogue
in Lemma \ref{trl} of the hypothesis that $G$ has few triangles, while the conclusion $\bigwedge_{e \in \ideal_\max} F_e = \emptyset$ is
the analogue of the conclusion that the modified graph $G'$ is triangle-free.  

\begin{proof}[of Lemma \ref{trl} assuming Theorem \ref{compact-triangle}]
Suppose for contradiction that Lemma \ref{trl} failed.  Then we can find an $0 < \eta \leq 1$ and 
sequence $n^{(m)}$ of integers, a sequence of graphs $G^{(m)} = (V^{(m)}, E^{(m)})$ with $|V^{(m)}|=n^{(m)}$, such that the $G^{(m)}$ have asymptotically vanishing number of triangles,
\begin{equation}\label{lim}
\lim_{m \to \infty} \frac{1}{(n^{(m)})^3} 
| \{ (x_1,x_2,x_3) \in (V^{(m)})^3: \{x_1,x_2\}, \{x_2,x_3\}, \{x_3,x_1\} \in E^{(m)} \} | = 0
\end{equation}
but such that each of the $G^{(m)}$ cannot be made triangle-free without deleting at least $\eta (n^{(m)})^2$ edges.  (One could make the decay
rate in \eqref{lim} more rapid, as in the proof of Theorem \ref{sz-quant}, but we will find no need to do so here.)  In particular, $G^{(m)}$ contains at least one triangle, and hence the expression inside the limit of \eqref{lim} is at least $1/(n^{(m)})^3$.  This implies that
\begin{equation}\label{nmgrow}
n^{(m)} \to \infty \hbox{ as } m \to \infty.
\end{equation}

Now let $(\Omega, \B_{\max})$ be the graph universal space introduced in Definition \ref{graph-universal}, with the attendant events $A_{n,m}$, regular algebra $\B_{\reg}$, factors $\B_I$, and $S_\infty$ group action. Let $\P^{(m)}$ be the probability measure
on $(\Omega \times \Omega^{(m)}, \B_{\max} \times \B_{\max}^{(m)})$ defined in Definition \ref{graph-embed}, and let $\P^{(\infty)}$ be a limiting
measure as constructed in the graph correspondence principle (Proposition \ref{gcp}).
From \eqref{lim} we have
$$ \lim_{m \to \infty} \P^{(m)}( A_{1,2} \wedge A_{2,3} \wedge A_{3,1} ) = 0$$
and hence by \eqref{weak-conv}
$$ \P^{(\infty)}( A_{1,2} \wedge A_{2,3} \wedge A_{3,1} ) = 0.$$
We will apply Theorem \ref{compact-triangle} on the universal space $(\Omega, \B_{\max}, \P^{(\infty)})$ 
with $J := \{1,2,3\}$, $\ideal_\max := \{ e: |e| \leq 2 \}$, $E_e$ set equal to $A_{i,j}$ if
$e = \{i,j\}$ for some $ij=12,23,31$, and $E_e = \Omega$ otherwise, and with $\B_e$ set equal to $\B_{e \cup \{4,5,\ldots\}}$ for all
$e \in \ideal_\max$.  The nesting and regularisability properties required for Theorem \ref{compact-triangle}
are obvious, while the independence properties follow from Lemma \ref{irl}.  We can thus invoke the theorem and
find regular events $F_e \in \B_{e \cup \{4,5,\ldots\}} \cap \B_\reg$ for $e \in \ideal_\max$ with
\begin{equation}\label{puff}
\P^{(\infty)}( E_e \backslash F_e ) < \eta/100 \hbox{ for } e \in \ideal_\max
\end{equation}
such that
\begin{equation}\label{fff}
\bigwedge_{e \in \ideal_\max} F_e = \emptyset.
\end{equation}

It is convenient to eliminate the lower order components $\emptyset, \{1\}, \{2\}, \{3\}$ of the ideal $\ideal_\max$.  For $ij=12,23,31$, define
$F'_{i,j} := F_{\{i,j\}} \wedge F_{\{i\}} \wedge F_{\{j\}} \wedge F_{\emptyset}$.  Then the $F'_{i,j}$ are regular, and (by monotonicity)
we have $F'_{i,j} \in \B_{\{i,j,4,5,\ldots\}}$.  From \eqref{puff} and the choice of the $E_e$ we have
\begin{equation}\label{sum}
 \P^{(\infty)}( A_{i,j} \backslash F'_{i,j} ) < \eta / 10 \hbox{ for } ij = 12, 23, 31
 \end{equation}
while from \eqref{fff} we have
\begin{equation}\label{fff2}
F'_{1,2} \wedge F'_{2,3} \wedge F'_{3,1} = \emptyset.
\end{equation}
Now we reinstate the ``hidden variables'' $x^{(m)}_1, x^{(m)}_2, \ldots$ by viewing $\P^{(m)}$ once again as a probability measure on
the product space $(\Omega \times \Omega^{(m)}, \B_{\max} \times \B_\max^{(m)})$; in particular $A_{i,j}$ is now identified with the event
that $\{x^{(m)}_i,x^{(m)}_j\}$ lies in the graph $G^{(m)}$.  Now because $A_{i,j}$ and $F'_{i,j}$ are regular, the quantity $\P^{(m)}( A_{i,j} \backslash F'_{i,j} )$ is the probability of an event involving only finitely many of the random vertices $x^{(m)}_i$ of $V^{(m)}$; let us say that it involves only the vertices $x^{(m)}_1, \ldots, x^{(m)}_N$ (note that $N$ will be independent of $m$, depending only on the complexity of the event $F'_{ij}$).  By increasing $N$ if necessary we may assume $N > 3$.  Recall that $F'_{ij}$ depends only on the 
vertices $x^{(m)}_i, x^{(m)}_j$ and $x^{(m)}_4,\ldots,x^{(m)}_N$.
For any fixed values of $x^{(m)}_4,\ldots,x^{(m)}_N$, let us say that a vertex pair $\{ x, y \} \subset  V^{(m)}$ is \emph{good} if for each $ij=12,23,31$, the event
$F'_{ij}$ holds true whenever $x,y$ are substituted for either $x^{(m)}_i, x^{(m)}_j$ or $x^{(m)}_j, x^{(m)}_i$.  
Now define the random subgraph $(G')^{(m)} =  (V^{(m)}, (E')^{(m)})$ of $G^{(m)} = (V^{(m)}, E^{(m)})$ by setting
$(E')^{(m)}$ to be all the good pairs $\{x,y\}$ in $E^{(m)}$; this graph depends on the random variables $x^{(m)}_4,\ldots,x^{(m)}_N$.  From \eqref{sum}
we see that
$$ \E |E^{(m)} \backslash (E')^{(m)}| < \eta (N^{(m)})^2.$$
Also, we observe that regardless of the values of $x^{(m)}_3,\ldots,x^{(m)}_N$, the graph $(G')^{(m)}$ almost surely cannot contain 
any triangles, as this would contradict \eqref{fff}.  But by the pigeonhole principle we can find a deterministic representative of the random 
graph  $(G')^{(m)}$ for which
$$ |E^{(m)} \backslash (E')^{(m)}| < \eta (N^{(m)})^2,$$
and so we have made $G^{(m)}$ triangle-free by removing fewer than $\eta (N^{(m)})^2$ edges, a contradiction that establishes Lemma \ref{trl}.
\end{proof}

\begin{remark} The same arguments in fact give a \emph{subgraph removal lemma}, in which the triangle is replaced by another fixed subgraph.  The proof is the same, it is only the downset $\ideal_\max$ (and some minor numerical factors in the argument) which change significantly.  But in all these cases, the elements in the downset will only have cardinality at most two.  We will not give the details here since they will be subsumed by the hypergraph removal lemma in Theorem \ref{hyper-remove}.
The higher order cases of Theorem \ref{compact-triangle}, involving sets $e$ of three or more elements, do not actually get used in graph theory (which is ultimately only concerned with finite boolean combinations of relations that involve at most 
two vertices at a time), and are only of importance for hypergraph theory (in which one must now consider combinations of relations, each of which involve three or more vertices).
\end{remark}

As observed in \cite{tao-hyper}, there is a slightly stronger version of the triangle removal lemma which gives some further complexity information on $G'$, at the expense of conceding that $G'$ need not be a subgraph of $G$.  More precisely, we have

\begin{lemma}[Strong triangle removal lemma]\label{trs} Let $G = (V,E)$ be an undirected graph with $|V|=n$ vertices.  Suppose that $G$ contains fewer than $\delta n^3$ triangles for some $0 < \delta \leq 1$.  Then one can find a triangle-free graph $G'$ with $G \backslash G'$ containing fewer than $o_{\delta \to 0}(n^2)$ edges.  Furthermore, there exists a partition of $V$ into $O_{\delta}(1)$ components, such that when restricted to the edges joining any two of these partitions (which could be equal), then $G'$ is either a complete graph or an empty graph.
\end{lemma}

This stronger version of the lemma is a by-product of the usual proof of Lemma \ref{trl}, as the graph $G'$ is constructed by excluding certain bad
pairs of Szemer\'edi cells from the graph $G$.  It turns out that the infinitary approach can also yield this stronger lemma without much difficulty.

\begin{proof}  We again argue by contradiction.  But this time, the contradiction hypothesis yields a more complicated statement.
More precisely, if Lemma \ref{trs} failed, then we can find an $0 < \eta \leq 1$ and 
sequence $n^{(m)}$ of integers, a sequence of graphs $G^{(m)} = (V^{(m)}, E^{(m)})$ with $|V^{(m)}|=n^{(m)}$, and a sequence $M^{(m)}$ tending to infinity as $m \to \infty$, such that \eqref{lim} holds, but such that there does not exist any triangle-free graph $G'$ for which 
$G \backslash G'$ has fewer than $\eta (n^{(m)})^2$ edges, \emph{and} for which there exists a partition of $V^{(m)}$ into $M^{(m)}$ or fewer 
components, such that when restricted to the edges joining any two of these cells, $G'$ is either the complete graph or the
empty graph.

We now repeat all the arguments used to prove Lemma \ref{trl}, until we get to the point where we have created regular events $F'_{i,j} \in
\B_{\{i,j,4,5,\ldots\}}$ obeying \eqref{sum} and \eqref{fff2}.  Now we insert an additional step to lower the complexity of the
events $F'_{i,j}$.  Observe that $\B_{\{i,j,4,5,\ldots\}}$ is generated by the factor $\B_{\{i,4,5,\ldots\}} \vee \B_{\{j,4,5,\ldots\}}$,
together with the additional event $A_{i,j}$.  Thus we can write $F'_{i,j} \wedge A_{i,j} = F''_{i,j} \wedge A_{i,j}$ for some
$F''_{i,j} \in \B_{\{i,4,5,\ldots\}} \vee \B_{\{j,4,5,\ldots\}}$.  From \eqref{sum} we have
$$   \P^{(\infty)}( A_{i,j} \backslash F''_{i,j} ) < \eta / 10 \hbox{ for } ij = 12, 23, 31.$$
Now we argue that we still have the analogue of \eqref{fff2}, namely
\begin{equation}\label{f3}
 F''_{1,2} \wedge F''_{2,3} \wedge F''_{3,1} = \emptyset.
 \end{equation}
From \eqref{fff2} we already have
\begin{equation}\label{f4}
 (F''_{1,2} \wedge F''_{2,3} \wedge F''_{3,1}) \cap (A_{1,2} \wedge A_{2,3} \wedge A_{3,1}) = \emptyset.
\end{equation}
But the regular event $F''_{1,2} \wedge F''_{2,3} \wedge F''_{3,1}$ is a boolean combination of finitely many events $A_{i,j}$, where at most
one of the $i,j$ lie in $\{1,2,3\}$.  In other words, this combination does not involve $A_{1,2} \wedge A_{2,3} \wedge A_{3,1}$.  Thus if \eqref{f3} failed, so that there was an infinite graph $(\N, E_\infty)$ lying in $ F''_{1,2} \wedge F''_{2,3} \wedge F''_{3,1} $, we could modify
the graph $E_\infty$ on the edges $\{1,2\}, \{2,3\}, \{3,1\}$ so that it also lies in the set in \eqref{f4}, a contradiction.  

To summarise, we can safely replace $F'_{i,j}$ by the lower complexity event $F''_{i,j}$.  
Now we continue the argument in the proof of Lemma
\ref{trl} with this replacement, but define the edges of $(G')^{(m)}$ to be all the good pairs $\{x,y\}$ in $V^{(m)}$ rather than in $E^{(m)}$.
This means that $(G')^{(m)}$ is no longer a subgraph of $(G)^{(m)}$, but the property of being good is determined entirely by the regular
events $F''_{i,j}$, which in turn only involve finitely many events $A_{i,j}$ with at most one of the $i,j$ lying in $\{1,2,3\}$.  Inspecting the definition of a good pair, we see that for fixed $x^{(m)}_4,\ldots,x^{(m)}_N$ for $N$ sufficiently large, the graph $(G')^{(m)}$ has bounded complexity, in the sense that there is a partition of $V^{(m)}$ into $M$ cells, for some $M$ depending only on the $F''_{i,j}$, such that when restricted to the edges joining any two of these cells, $G'$ is either the complete graph or the
empty graph.  But for $m$ sufficiently large we have $M^{(m)} > M$, and so we attain the same contradiction as before.
\end{proof}

\section{The uniform intersection property}\label{uip-sec}

We now build the machinery necessary to prove the infinitary hypergraph removal lemma (Theorem \ref{compact-triangle}).  Again,
we will be motivated by the example from ergodic theory.  In Furstenberg's proof \cite{furst}, \cite{furstenberg}, \cite{furst-book} of
the Furstenberg recurrence theorem (Theorem \ref{furst-thm}), the proof proceeded by a kind of induction on factors.  Let us say that a factor
$\B$ of $\B_{\max}$ obeys the \emph{uniform multiple recurrence (UMR) property} if the conclusion of Theorem \ref{furst-thm} holds true whenever $E \in \B$ and $\P(E) > 0$.  Thus for instance the trivial factor $\{ \emptyset, \Omega\}$ has the UMR property.  One then shows that the UMR property is preserved under three operations: weakly mixing extensions; limits of chains; and compact (or finite rank) extensions. An application of Zorn's lemma\footnote{Actually, to establish Theorem \ref{furst-thm} for a fixed $k$, one only needs to apply the limits-of-chains step a finite number of times depending on $k$, at which point one reaches a factor which is \emph{characteristic} for the maximal factor $\B_{\max}$, at which point one can jump directly to $\B_{\max}$ without using Zorn's lemma.  Thus the proof of the Furstenberg recurrence theorem does not actually require the axiom of choice, and indeed the original proof in \cite{furst} did not use this axiom.} then allows one to conclude that the maximal factor $\B_{\max}$ obeys the UMR property.

We will adopt a similar strategy here, based around a certain property of families of factors which we call the 
\emph{uniform intersection property (UIP)}.  This property is again trivial for very small families, and will be preserved under the same three operations of weakly mixing extensions, limits of chains, and finite rank extensions.  Because of the finiteness of $J$ in Theorem \ref{compact-triangle}, we will only need to apply these operations finitely often, and will not require Zorn's lemma.  However it does seem likely that there are extensions of this theorem to the case when $J$ is infinite, and (more interestingly) to the case where the sets $e$ in $\ideal_\max$ can be unbounded or even countably infinite.  We will not pursue this matter here.

We begin by stating the UIP.

\begin{definition}[Uniform intersection property]  Let $(\Omega, \B_{\max}, \P)$ be a probability space, and let $\B_{\reg}$ be an algebra in $\B_{\max}$.
We say that a tuple $(\B_i)_{i \in I}$ of factors has the \emph{uniform intersection property} (UIP) if the following holds: given any tuple $(E_i)_{i \in I}$ of events $E_i \in \B_i$ with $\P( \bigwedge_{i \in I} E_i ) = 0$, and given any $\eps > 0$, there exists a tuple $(F_i)_{i \in I}$ of regular events $F_i \in \B_i \wedge \B_\reg$ with $\P(E_i \backslash F_i) \leq \eps$ for each $i \in I$ such that $\bigwedge_{i \in I} F_i = \emptyset$.
\end{definition}

\begin{remark} Roughly speaking, the UIP asserts that if events $E_i$ from separate factors $\B_i$ have a null intersection, then this fact can be almost entirely ``explained'' by \emph{regular} events $F_i \in \B_i$ which have \emph{empty} intersection. Thus, for instance, the conclusion of Theorem \ref{compact-triangle} is simply that the tuple $(\B_I)_{I \in \ideal_\max}$ obeys the UIP.
\end{remark}

Before continuing, let us illustrate the UIP with a few simple examples.  All of these examples take place in some probability space $(\Omega, \B_{\max}, \P)$ with an algebra $\B_\reg$ of regular events.  We say that a factor $\B$ is \emph{regularisable} if it can be generated by at most countably many regular events.

\begin{example}\label{empty} The empty tuple $()$ obeys the UIP in a vacuous sense (the hypothesis $\P( \bigwedge_{i \in I} E_i ) = 0$ is impossible to satisfy). 
\end{example}

\begin{example}\label{singleton} Let $\B$ be a factor.  Then 
the singleton tuple $( \B )$ trivially has the UIP (indeed one can even take $\eps=0$ and $F_\B = \emptyset$ in this case).
\end{example}

\begin{example} Let $\B$ be a regularisable factor.  Then the $2$-tuple $(\B,\B)$ has the UIP.  Indeed, if $E, E' \in \B$
were such that $\P( E \wedge E' ) = 0$, then from Lemma \ref{product-approx} one can find regular events $\tilde F, \tilde F' \in \B$ which are
$\eps/3$-close to $E$, $E'$ respectively.  By the triangle inequality this implies that $\P( \tilde F \wedge \tilde F' ) \leq 2\eps/3$.
If we set $F := \tilde F \backslash \tilde F'$ and $F' := \tilde F' \backslash \tilde F$ then we see that $F,F'$ are regular events in 
$\B$ with $F \wedge F' = \emptyset$, while from the triangle inequality $\P( E \backslash F ), \P( E' \backslash F' ) \leq \eps$, and the claim
follows.  For a generalization of this argument, see Lemma \ref{repeat} below.
\end{example}

\begin{example}\label{trivia} The trivial factor $\{\emptyset, \Omega\}$ has no impact on the UIP. More precisely,
a tuple $(\B_i)_{i \in I}$ obeys the UIP if and only if $(\B_i)_{i \in I} \uplus (\{\emptyset,\Omega\})$ also obeys the UIP, where we use $\uplus$ to denote the concatenation of tuples.
\end{example}

\begin{example}\label{independent} Let $( \B_1, \ldots, \B_l )$ be a tuple of jointly independent factors.  
Then $( \B_1, \ldots, \B_l )$ has the UIP.  Indeed, if $E_j \in \B_j$
then by joint independence we have $\P( \bigwedge_{1 \leq j \leq l} E_j ) = \prod_{j=1}^l \P( E_j )$.  Thus if $\bigwedge_{1 \leq j \leq l} E_j$
is a null event, then one of the $E_j$, say $E_{j_0}$, must also be a null event.  The claim then follows by letting $F_{j_0} = \emptyset$
and letting all the other $F_j$ be the full event $\Omega$.  For a more sophisticated version of this argument, see Lemma \ref{weakmix} below.
\end{example}

\begin{example}  Let $\Omega$ be the unit interval $[0,1]$ with Lebesgue measure, let $\B_\reg$ consist of all the finite unions of intervals (open, closed, or half-open), let $\B_1$ be the factor generated by the event $E_1 := [0,1/2]$,
and let $\B_2$ be the factor generated by the event $E_2 := [1/2,1]$.  Then $(\B_1,\B_2)$ does not have the UIP.
However if one modifies $\B_2$ to be the factor generated by $(1/2,1]$ instead of $[1/2,1]$, then the UIP is restored.  Thus the UIP is sensitive to modification of the underlying factors by null events.  (On the other hand, the events $E_i$ themselves can be modified by null events within $\B_i$ without any impact to the UIP.)
\end{example}

\begin{example} Let $\B_1, \B_2$ be finite $\sigma$-algebras, and let $\B$ be another $\sigma$-algebra, such that $(\B_1, \B_2, \B)$ has the UIP.  Then $(\B_1 \vee \B_2, \B)$ also has the UIP.  To see this, let $E_{12} \in \B_1 \vee \B_2$ and $E \in \B$ be such that $\P( E_{12} \vee \B ) = 0$.
Since $\B_1$ and $\B_2$ are finite, we can write $E_{12}$ as the union of $M$ events of the form $E_{1,m} \wedge E_{2,m}$ for $1 \leq m \leq M$
for some finite $M$ and some events $E_{1,m} \in \B_1, E_{2,m} \in \B_2$.  From the UIP hypothesis we can find regular events $F_{1,m} \in \B_1$,
$F_{2,m} \in \B_2$, $F_m \in \B$ with $F_{1,m} \wedge F_{2,m} \wedge F_m = \emptyset$ and $\P( E_{1,m} \backslash F_{1,m} ), \P( E_{2,m} \backslash F_{2,m}), \P( E \backslash F_m ) \leq \eps/M$.  If we then set $F_{12} := \bigvee_{m=1}^M (F_{1,m} \wedge F_{2,m})$ and $F := \bigwedge_{m=1}^M F_m$
then the claim follows.  For a generalization of this argument, see Lemma \ref{finite} below.
\end{example}

\begin{remark}\label{cantelli}  If $(\B_i)_{i \in I}$ has the UIP, then given any tuple $(E_i)_{i \in I}$ of events $E_i \in \B_i$ such that $\bigwedge_{i \in I} E_i$ is a null event, then there exists a tuple $(G_i)_{i \in I}$ of null events $G_i \in \B_i$ which cover the 
null event $\bigwedge_{i \in I} E_i$.  This
follows from applying the UIP with $\eps = 2^{-n}$ (say) to obtain events $(F_{i,n})_{i \in I}$ with $\P( E_i \backslash F_{i,n} ) \leq 2^{-n}$ 
for all $i \in I$ and
$\bigwedge_{i \in I} F_{i,n} = \emptyset$, and
then letting $G_i$ be the event $E_i$ holds, but that $F_{i,n}$ fails for infinitely many $n \geq 1$; the claim $\bigwedge_{i \in I} E_i \subseteq \bigvee_{i \in I} G_i$ then
follows from the pigeonhole principle, while the claim that $G_i$ is null follows from the Borel-Cantelli lemma.  Unfortunately the $G_i$ will in general not be regular, and so this consequence of the UIP, while simple to state, is not useful for applications.
\end{remark}

We now develop the general tools that we shall use to deduce the UIP for complex tuples from the UIP for simpler tuples.
We first show that repetitions do not affect the UIP so long as the $\sigma$-algebra being repeated is regularisable.  We use $(\B_i)_{i \in I}\uplus (\B_j)_{j \in J} = (\B_i)_{i \in I \uplus J}$ to denote the concatentation of two tuples, where $I\uplus J$ is the disjoint union of $I$ and $J$ (thus one may have to relabel the index set of $I$ or $J$ in order to define this concatenation).

\begin{lemma}\label{repeat}  Let $(\B_i)_{i \in I}$ be a tuple of $\sigma$-algebras, and let $\B$ be a regularisable $\sigma$-algebra.  Then
$(\B_i)_{i \in I} \uplus (\B)$ has the UIP if and only if $(\B_i)_{i \in I} \uplus (\B,\B)$.
\end{lemma}

\begin{proof} First suppose that $(\B_i)_{i \in I} \uplus (\B,\B)$ has the UIP.  Then if $E_i \in \B_i$ and $E \in \B$ are such that
$\P( E \wedge \bigwedge_{i \in I} E_i ) = 0$, then by the UIP hypothesis (inserting a dummy event $\Omega$ for the extra copy of $\B$)
we can find regular events $F_i \in \B_i$ and $F', F'' \in \B$ such that $F' \wedge F'' \wedge \bigwedge_{i \in I} F_i = \emptyset$ and
$$\P( E_i \backslash F_i ), \P( E \backslash F' ), \P( \overline{F''} ) \leq \eps/2 \hbox{ for all } i \in I.$$
The claim then follows by setting $F := F' \wedge F''$.

Now suppose conversely that $(\B_i)_{i \in I} \uplus (\B)$ has the UIP.  Then if $E_i \in \B_i$ and $E, E' \in \B$ are such that
$\P( E \wedge E' \wedge \bigwedge_{i \in I} E_i ) = 0$ then by the UIP hypothesis (replacing $E$ and $E'$ by the single event $E \wedge E'$)
one can find regular events $F_i \in \B_i$ and $\tilde F \in \B$ such that $\tilde F \wedge \bigwedge_{i \in I} F_i = \emptyset$ and
$$\P( E_i \backslash F_i ), \P( (E \wedge E') \backslash \tilde F ) \leq \eps/3 \hbox{ for all } i \in I.$$
Now since $\B$ is regularisable, we see from Lemma \ref{product-approx} that every event in $\B$ is $\eps$-close to a regular event in $\B$ for any
$\eps > 0$.  In particular, we can find regular events $\tilde E, \tilde E' \in \B$ which are $\eps/3$-close to $E$ and $E'$ respectively.  If one
then sets $F := (\tilde E \backslash \tilde E') \vee \tilde F$ and $F' := (\tilde E' \backslash \tilde E) \vee \tilde F$, we see from
the triangle inequality that $\P( E \backslash F ), \P( E' \backslash F' ) \leq \eps$, and that $F \wedge F' \wedge \bigwedge_{i \in I} F_i = \emptyset$, and the claim follows.
\end{proof}

Now we give the three major extendability properties of the UIP, under weakly mixing extensions, finite rank extensions, and limits of chains.
We begin with the analogue of the weakly mixing extension property, which says that one can extend any member of a tuple without destroying the 
UIP, as long as the extension is relatively independent of all the other factors in the tuple.

\begin{lemma}[Weakly mixing extensions]\label{weakmix}  Let $(\B_i)_{i \in I}$ be a tuple of $\sigma$-algebras, and let $\B$ be an additional $\sigma$-algebra such that $(\B_i)_{i \in I} \uplus (\B)$ obeys the UIP.  Let $\B'$ be a extension of $\B$ which is relatively independent of $\bigvee_{i \in I} \B_i$ over $\B$.  Then $(\B_i)_{i \in I} \uplus (\B')$ also obeys the UIP.
\end{lemma}

\begin{proof}  Let $E_i \in \B_i$ for $i \in I$ and $E' \in \B'$ be events such that $\P( E' \wedge \bigwedge_{i \in I} E_i ) = 0$.  
We rewrite this as
$$ \E( \I(E') \prod_{i \in I} \I(E_i) ) = 0.$$
The first factor is measurable in $\B'$, while the second factor is measurable in $\bigvee_{i \in I} \B_i$.  Since these two $\sigma$-algebras are relatively independent over $\B$, we 
may use \eqref{eff} and conclude
that $\P(E'|\B) \prod_{i \in I} \I(E_i) = 0$ almost surely.
Let $E \in \B$ be the support of $\P(E'|\B)$ (which is determined only up to a null event in $\B$), then
$E \wedge \bigwedge_{i \in I} E_i$ is a null event.  
Applying the UIP hypothesis, we can find regular events $F_i \in \B_i$ for $i \in I$ and a regular
event $F \in \B$ such that $\P(E_i \backslash F_i), \P(E \backslash F) \leq \eps$ and
$$ F \vee \bigvee_{i \in I} F_i = \emptyset.$$
We then set $F' := F$.  We will be done as soon as we check that $\P( E' \backslash F' ) \leq \eps$,
which will follow if we can show that $\P( E' \backslash E ) = 0$.  But
$$ \P( E' \backslash E ) = \P( \I( E') (1 - \I( E )) ) =
\P( \P( E' | \B ) (1 - \I( E )) )$$
since $1 - \I(E)$ is $\B$-measurable.  But this vanishes by the definition of $E'$. 
\end{proof}

Next we turn to the preservation of the UIP under compact extensions (or more accurately ``finite rank extensions''), which assert that one can extend any given element of a tuple by finite factors of other elements in the tuple (destroying those elements in the process). 

\begin{lemma}[Finite rank extensions]\label{finite}  Let $\B_0$ be a $\sigma$-algebra, let $\B_1,\ldots,\B_l$ be factors of $\B_0$, and let $\B'_1,\ldots,\B'_l$ be finite $\sigma$-algebras for some $l \geq 1$.  Let $(\tilde \B_i)_{i \in I}$ be an additional tuple of $\sigma$-algebras.  Then if $(\tilde \B_i)_{i \in I} \uplus ( \B_1 \vee \B'_1, \ldots, \B_l \vee \B'_l, \B_0 )$ has the UIP, then
$(\tilde \B_i)_{i \in I} \uplus ( \B_0 \vee \B'_1 \vee \ldots \vee \B'_l )$ also has the UIP.
\end{lemma}

\begin{proof}  Write $\B_* := \B_0 \vee \B'_1 \vee \ldots \vee \B'_l$.  
Let $E_*$ be an event in $\B_*$, and let $\tilde E_i \in \tilde \B_i$ for $i \in I$, be such that $\P( E_* \wedge \bigwedge_{i \in I} \tilde E_i ) = 0$.
Since $\B'_1,\ldots,\B'_l$ are finite, we can write $E_*$ as the finite union of events
$$ E_* = E_{*,1} \vee \ldots \vee E_{*,M}$$
for some $M \geq 1$, where each $E_{*,m}$ has the form
$$ E_{*,m} = E_{0,m} \wedge E_{1,m} \wedge \ldots \wedge E_{l,m}$$
for some events $E_{0,m} \in \B_0$ and $E_{j,m} \in \B'_j$ for $1 \leq j \leq l$.  For each $1 \leq m \leq M$, we have
$\P( E_{*,m} \wedge \bigwedge_{i \in I} \tilde E_i ) = 0$ and hence
$$ \P( E_{0,m} \wedge E_{1,m} \wedge \ldots \wedge E_{l,m} \wedge \bigwedge_{i \in I} \tilde E_i ) = 0.$$
Observe that $E_{j,m} \in \B_j \vee \B'_j$ for $1 \leq j \leq l$, and hence by the UIP hypothesis we may find regular events
$F_{0,m} \in \B_0$, $F_{j,m} \in \B_j \vee \B'_j \subseteq \B_*$ for $1 \leq j \leq l$, and $\tilde F_{i,m} \in \tilde \B_i$ for $i \in I$ such that
$$ \P( E_{0,m} \backslash F_{0,m} ), \P( E_{j,m} \backslash F_{j,m} ), \P( \tilde E_i \backslash \tilde F_{i,m} ) \leq \frac{\eps}{M(l+1)}$$
for $1 \leq j \leq l$ and $i \in I$, and
$$ \bigwedge_{j=1}^l F_{j,m} \wedge F_{0,m} \wedge \bigwedge_{i \in I} \tilde F_{i,m} = \emptyset.$$
Thus if we set 
$F_* := \bigvee_{m=1}^M (\bigwedge_{j=1}^l F_{j,m} \wedge F_{0,m})$ and $\tilde F_i := \bigwedge_{m=1}^M \tilde F_{i,m}$ for $i \in I$
then $F_* \in \B_*$ and $\tilde F_i \in \tilde \B_i$ are regular events, and $F_* \wedge \bigwedge_{i \in I} \tilde F_i = \emptyset$, and
$$ \P( E_* \backslash F_*) \leq \sum_{m=1}^M \P( E_{*,m} \backslash (\bigwedge_{j=1}^l F_{j,m} \wedge F_{0,m}) ) 
\leq M (l+1) \frac{\eps}{M(l+1)} = \eps$$
and
$$ \P( \tilde E_i \backslash \tilde F_i ) \leq \sum_{m=1}^M \P( \tilde E_i \backslash \tilde F_{i,m} ) \leq M \frac{\eps}{M(l+1)} \leq \eps$$
and the claim follows.
\end{proof}

Finally, we consider the preservation of the UIP under limits of chains assuming a certain relative independence property.

\begin{lemma}[Limits of chains]\label{zorn} Let $A$ be a totally ordered set, let $I$ be a finite index set, and for each $\alpha \in A$ let $(\B_{\alpha,i})_{i \in I}$ be a tuple of $\sigma$-algebras obeying the UIP, which is increasing in the sense that $\B_{\alpha,i}$ is a factor of $\B_{\beta,i}$ whenever $\alpha < \beta$ and $i \in I$.  Let $\B_i := \bigvee_{\alpha \in A} \B_{\alpha,i}$,
and suppose that whenever $i \in I$ and $\alpha \in A$,
the $\sigma$-algebras $\B_i$ and $\bigvee_{j \in I \backslash \{i\}} \B_{\alpha,j}$ are relatively independent over $\B_{\alpha,i}$.
Then the tuple $(\B_i)_{i \in I}$ also obeys the UIP.
\end{lemma}

\begin{proof}  
Let $E_i \in \B_i$ for $i \in I$ be such that
$\P( \bigwedge_{i \in I} E_i ) = 0$.  From Corollary \ref{product-monotone}, we see that
for each $i \in I$ there exists an $\alpha \in A$ such that $E_i$ is $\eps/4(|I|+1)^2$-close to an event in $\B_{\alpha,i}$.  
Since there are only finitely many $i$, we can make this $\alpha$ uniform in $i$.  This
implies in particular that $\| \I( E_i ) - \P( E_i | \B_{\alpha,i} ) \|_{L^2} \leq \eps^{1/2}/2(|I|+1)$ for all $i \in I$,
since the orthogonal projection $\P( E_i  | \B_{\alpha,i} )$ is the nearest $\B_{\alpha,i}$-measurable random variable to $\I( E_i )$ 
in the $L^2$ metric. 

Let $E_{\alpha,i} \in \B_{\alpha,i}$ denote the event that $\P( E_i  | \B_{\alpha,i} ) > \frac{|I|}{|I|+1}$ (this event is only defined
up to null events in $\B_{\alpha,i}$).   Then by Chebyshev's inequality we have
\begin{align*}
\P( E_i \backslash E_{\alpha,i} ) &\leq \P\left( |\I( E_i ) - \P( E_i | \B_{\alpha,i} )| > \frac{1}{|I|+1} \right) \\
&\leq \eps/2
\end{align*}
for each $i \in I$.  Now let $A$ denote the event $A := \bigwedge_{i \in I} E_{\alpha,i}$.  Since $\P(\bigwedge_{i \in I} E_i) = 0$, we have
\begin{equation}\label{omega-climb}
 \P(A) \leq \sum_{i \in I} \P( A \backslash E_i ).
 \end{equation}
On the other hand, we have
$$ \P( A \backslash E_i) = \E( \I(E_{\alpha,i} \backslash E_i) \prod_{j \in I \backslash \{i\}} \I( E_{\alpha,j} ) ).$$
Using the relative independence hypothesis, we conclude
$$ \P( A \backslash E_i) = \E( \P(E_{\alpha,i} \backslash E_i | \B_{\alpha,i}) \prod_{j \in I \backslash \{i\}} \I( E_{\alpha,j} ) ).$$
But by definition of $E_{\alpha,i}$ we have
$$ \P( E_{\alpha,i} \backslash E_i | \B_{\alpha,i} ) = \I( E_{\alpha,i} ) (1 - \P( E_i | \B_{\alpha,i} ) ) \leq \frac{1}{|I|+1} \I( E_{\alpha,i} ) $$
and hence
$$ \P( A \backslash E_i ) \leq \frac{1}{|I|+1} \E( \I(E_{\alpha,i} ) \prod_{j \in I \backslash \{i\}} \I( E_{\alpha,j} ) )
= \frac{1}{k+1} \P(A).$$
Inserting this back into \eqref{omega-climb} we conclude that $\P(A) \leq \frac{|I|}{|I|+1} \P(A)$, and hence that $A$ is a null event.
By definition of $A$ and the UIP hypothesis, we may thus find regular events $F_{\alpha,i} \in \B_{\alpha,i}$ for all $i \in I$
with $\P(E_{\alpha,i} \backslash F_{\alpha,i}) \leq \eps/2$ such that $\bigwedge_{i \in I} F_{\alpha,i} = \emptyset$.
By the triangle inequality we thus have $\P(E_i \backslash F_{\alpha,i}) \leq \eps/2$.  The claim now follows by setting
$F_i := F_{\alpha,i}$.
\end{proof}

\section{Proof of the infinitary hypergraph removal lemma}\label{hypersec}

We are now ready to prove the hypergraph removal lemma.  Fix the probability space $(\Omega, \B_{\max},\P)$, the algebra $\B_\reg$ of regular events,
the finite set, $J$, the downset $\ideal_\max$, and the factors $\B_I$ for $I \in \ideal_\max$ obeying the hypotheses in Theorem \ref{compact-triangle}.  For any sub-ideal $\ideal$ of $\ideal_\max$, let $\B(\ideal)$ denote the factor $\B(\ideal) := \bigvee_{I \in \ideal} \B_I$, thus $\B(\ideal)$ is a regularisable factor.  For any $I \in \ideal_\max$, define the \emph{principal ideal} $\langle I \rangle := \{ I': I' \subseteq I \}$; from the nesting property we see that $\B(\langle I \rangle) = \B_I$ for all $I \in \ideal_\max$. Thus our task is to show that
the tuple $(\B( \langle I \rangle ) )_{I \in \ideal_{\max}}$ obeys the UIP.  For inductive purposes, we will derive this claim from the following more general statement.  For any downset $\ideal$, we define the \emph{height} $\height(\ideal)$ of $\ideal$ to be the quantity
$\height(\ideal) := \sup \{ |e| : e \in \ideal \}$, with the convention that the empty ideal has height $-\infty$.

\begin{proposition}\label{propinduct} Let the hypotheses and notation be as above.  Let $d \geq 0$, and let $(\ideal_i)_{i \in I}$ be any finite tuple of sub-ideals of $\ideal_\max$ (possibly with repetitions), such that every ideal $\ideal_i$ has height at most $d$.  Then the tuple $(\B(\ideal_i))_{i \in I}$ obeys the UIP.
\end{proposition}

By taking $d$ sufficiently large (e.g. $d=|J|$) we obtain Theorem \ref{compact-triangle}.

\begin{proof}
We will prove Proposition \ref{propinduct} by an induction on $d$.  First consider the base case $d=0$.  Then the only ideals available are
the empty ideal $\{\}$, and the singleton ideal $\{\emptyset\}$; these correspond to the trivial factor $\{\emptyset, \Omega\}$ and the
regularisable factor $\B_\emptyset$.  The claim now follows from Examples \ref{empty}, \ref{singleton}, \ref{trivia} and Lemma \ref{repeat}.

Now suppose that $d \geq 1$, and that Proposition \ref{propinduct} has already been proven for $d-1$.  
First observe from Lemma \ref{repeat} that we may remove duplicates and assume that all the ideals $\ideal_i$ are distinct.

Given any $e \in \ideal_\max$ with $|e|=d$, we know that $\B_e$ is regularisable, hence we may write $\B_e = \bigvee_{n \geq 1} \B_{e,n}$
for some increasing sequence $\B_{e,1}\subseteq \B_{e,2} \subseteq \ldots$ of regularisable finite $\sigma$-algebras.  In particular,
we have $\B(\ideal_i) = \bigvee_{n \geq 1} \B_n(\ideal_i)$ for all $i \in I$, where
$$ \B_n(\ideal) := \bigvee_{e \in \ideal: |e|=d} \B_{e,n} \vee \B(\overline{\ideal})$$
and $\overline{\ideal}$ is the downset $\overline{\ideal} := \{ e \in \ideal: |e'| < d \}$; note that this ideal has height strictly 
less than $d$.

We need some relative independence properties of the factors $\B(\ideal)$.  We begin with

\begin{lemma}\label{crop} Let $\ideal, \ideal'$ be sub-ideals of $\ideal_\max$ height at most $d$ which do not have any common elements of order exactly $d$.  Then
$\B(\ideal)$ and $\B(\ideal')$ are relatively independent over $\B(\overline{\ideal})$.
\end{lemma}

\begin{proof} We will induct on the quantity $m := |\{ e \in \ideal: |e| = d \}|$, the number of top-order elements in $\ideal$.
If $m=0$ then $\ideal = \overline{\ideal}$ and the claim follows.  Now suppose that $m \geq 1$ and the
claim has already been established for $m-1$.  Let $e_d$ be an element of $\ideal$ with $|e_d|=d$, and let $\tilde \ideal:= \ideal \backslash \{e_d\}$.  From the induction hypothesis we already know that
$\B(\tilde \ideal)$ and $\B(\ideal')$
are relatively independent over $\B(\overline{\ideal})$.  Also, from the UIP hypothesis we know that
$\B(\langle e_d \rangle)$ and $\B(\tilde \ideal) \vee \B(\ideal')$
are relatively independent over $\B(\overline{\ideal})$.  Applying the gluing property (Proposition \ref{stability}(i)) we conclude
that the factors $\B(\tilde \ideal) \vee \B(\langle e_d \rangle)$ and $\B(\ideal')$
are relatively independent over $\B(\overline{\ideal})$.  Since the former factor is nothing more than
$\B(\ideal)$, the claim follows.
\end{proof}

As a consequence, we have

\begin{lemma} Let $i \in I$ and $n \geq 1$.  Then $\B(\ideal_i)$ and $\bigvee_{j \in I: j \neq i} \B_n(\ideal_j)$ are relatively
independent over $\B_n(\ideal_i)$.
\end{lemma}

\begin{proof} Observe that $\bigvee_{j \in I: j \neq i} \B_n(\ideal_j)$ is a factor of $\B_n(\ideal_i) \vee \B(\ideal')$, where
$\ideal'$ is the downset
$\ideal' := (\bigcup_{j \in I} \ideal_j) \backslash \{ e \in \ideal_i: |e|=d \}$.  Thus by monotonicity and absorption (Proposition \ref{stability}(i),
(ii)) it suffices to show that $\B(\ideal_i)$ and $\B(\ideal')$ are relatively independent over $\B_n(\ideal_i)$.  Since
factors do not affect relative independence (Proposition \ref{stability}(iv)), it suffices to show that $\B(\ideal_i)$ and $\B(\ideal')$ are relatively independent over
$\B(\overline{\ideal_i})$.  But this follows from Lemma \ref{crop}.
\end{proof}

From the above lemma and Lemma \ref{zorn}, we see that to close the induction hypothesis it suffices
to show that $(\B_n(\ideal_i))_{i \in I}$ obeys the UIP for all $n \geq 1$.

Let $k$ denote the number of ideals $\ideal_i$ whose height is exactly $d$.  First suppose that all the ideals $\ideal_i$ have height strictly
less than $d$.  Then $\B_n(\ideal_i) = \B(\ideal_i)$, and the claim follows from the induction hypothesis.

Now suppose that all the ideals $\ideal_i$ either have height strictly less than $d$, or are principal ideals (this is a ``weakly mixing'' case).  
We induct on the number of principal ideals of height $d$.  If there are no such ideals, then we are done by the preceding paragraph.  Since we have removed duplicates, we know that no two principal ideals present have any common elements of top order $d$.  Thus if $\ideal_i$ is a principal ideal, then $\B_n(\ideal_i)$ is relatively independent of $\bigvee_{j \in I: j \neq i} \B_n(\ideal_j)$ over $\B(\overline{\ideal_i})$.  Applying
Lemma \ref{weakmix}, it suffices for the purposes of checking the UIP to replace $\B_n(\ideal_i)$ with $\B_n(\overline{\ideal_i})$.  But this
follows from the (inner) induction hypothesis.

Finally, we consider the general case.  Let $k$ denote the number of ideals $\ideal_i$ of height $d$ which are not principal.  We have already dealt
with the case $k=0$, so suppose inductively that $k \geq 1$ and the claim has already been proven for $k-1$.  Let
Let $\ideal_{i_0}$ be an ideal of height $d$ which is not principal, and let $e_1,\ldots,e_l$ be the elements of $\ideal_{i_0}$ of
order $d$.  We can then split
$$ \B_n(\ideal_{i_0}) = \B_n( \langle e_1 \rangle ) \vee \ldots \vee \B_n( \langle e_l \rangle ) \vee \B(\overline{\ideal_i}).$$
Also observe that for $1 \leq j \leq l$ we have
$\B_n( \langle e_j \rangle ) = \B_{e_j,n} \vee \B( \overline{\langle e_j \rangle})$, and that $\B( \overline{\langle e_j \rangle})$ is a factor of $\B(\overline{\ideal_i})$.  Thus we may apply Lemma \ref{finite} and conclude
that in order to prove the UIP for $(\B_n(\ideal_i))_{i \in I}$, it suffices to do so for the tuple
$$ (\B_n(\ideal_i))_{i \in I \backslash \{i_0\}} \uplus ( 
\B_n( \langle e_1 \rangle ) , \ldots, \B_n( \langle e_l \rangle ), \B(\overline{\ideal_i})).$$
This tuple has one fewer non-principal degree $d$ ideal than the original tuple, and so the claim now follows from the (inner) induction hypothesis.
\end{proof}

\section{A hypergraph correspondence principle}\label{hyperspace}

We now generalise the graph correspondence principle developed in Section \ref{graphsec} to hypergraphs. To keep the exposition somewhat simple we shall restrict our attention to the principle for a single $d$-uniform hypergraphs, although there would be no difficulty extending this principle to systems of hypergraphs of varying uniformities and partite-ness.  The material here will be extremely analogous to Section \ref{graphsec}. Indeed, we could have deleted that section as being redundant, but we believe for pedagogical purposes that it is better to start with graphs before moving on to hypergraphs.

\begin{definition}[Hypergraphs]  Let $d \geq 0$.  If $V$ is a set, we let $\binom{V}{d} := \{ e \subset V: |e| = d \}$ denote the $d$-element subsets of $V$.  A \emph{$d$-uniform hypergraph} is a pair $G = (V, E)$, where $V$ is a non-empty set and $E \subset \binom{V}{d}$.  
\end{definition}

Note that a $2$-uniform hypergraph is the same concept as an undirected graph.  We will fix $d \geq 2$, and
consider the question of extracting an infinitary limit from a sequence $G^{(m)} = (V^{(m)}, E^{(m)})$ of $d$-uniform hypergraphs.  As before, we shall
need a universal space, an embedding into that space, and a correspondence principle.  We begin with the universal space.

\begin{definition}[Hypergraph universal space]\label{hypergraph-universal}  Fix $d \geq 2$.
Let $\Omega := 2^{\binom{\N}{d}} = \{ (\N, E_\infty): E_\infty \subset \binom{\N}{d} \}$ denote the space of all infinite $d$-uniform hypergraphs $(\N, E_\infty)$ on the natural numbers.  On this space $\Omega$, we introduce the events $A_e$ for all $e \in \binom{\N}{d}$ by
$A_e := \{ (\N,E_\infty) \in \Omega: e \in E_\infty \}$, and let $\B_{\max}$ be the $\sigma$-algebra generated by the $A_e$.   
We also introduce the regular algebra $\B_{\reg}$ generated by the $A_e$, thus these are the events that depend only only finitely many of 
the $A_e$.
For any $\sigma \in S_\infty$, we define the associated action on $\B_{\max}$ by mapping $\sigma: A_e \mapsto A_{\sigma(e)}$ and extending this to a $\sigma$-algebra isomorphism in the unique manner.
For any (possibly infinite) subset $I$ of $\N$, we define $\B_I$ to be the factor of $\B_{\max}$ generated by
the events $A_e$ for $e \in \binom{I}{d}$.
\end{definition}

Next, we need a way to embed every finite hypergraph into the universal space.

\begin{definition}[Hypergraph universal embedding]\label{hypergraph-embed}  Fix $d \geq 2$.  Let $m \geq 1$, and let $G^{(m)} = (V^{(m)}, E^{(m)})$ be a finite $d$-uniform hypergraph.
Let $(\Omega^{(m)}, \B_{\max}^{(m)}, \P^{(m)})$ be the probability space corresponding to the sampling of a countable sequence of iid random variables $x^{(m)}_1, x^{(m)}_2, \ldots \in V^{(m)}$ sampled independently and uniformly at random.  To every sequence
$(x^{(m)}_1,x^{(m)}_2, \ldots ) \in \Omega^{(m)}$ we associate an infinite $d$-uniform hypergraph $G^{(m)}_\infty = (\N, E^{(m)}_\infty) \in \Omega$ by setting
$$ E^{(m)}_\infty := \{ e \in \binom{\N}{d}: \{x^{(m)}_i: i \in e\} \in E^{(m)} \}.$$
This mapping from $\Omega^{(m)}$ to $\Omega$ is clearly measurable, since the inverse images of the generating events $A_e \in (\Omega, \B_{\max})$ 
are the events that $\{x^{(m)}_i: i \in e\}$ lie in $G^{(m)}$, which are certainly measurable in $\B^{(m)}_\max$.  
This allows us to extend the probability measure
$\P^{(m)}$ from $(\B_{\max}^{(m)}, \Omega^{(m)})$ to the product space $(\B_{\max} \times \B_{\max}^{(m)}, \Omega \times \Omega^{(m)})$ 
in a canonical manner, identifying the events $A_{e}$ with the events $\{ x^{(m)}_i: i \in e \} \in E^{(m)}$.  We shall abuse notation and
refer to the extended measure also as $\P^{(m)}$.
\end{definition}

As before we can verify the permutation invariance \eqref{perm}.  
By repeating the proof of the graph correspondence principle (Proposition \ref{gcp}) almost word-for-word, we obtain its counterpart for hypergraphs:

\begin{theorem}[Hypergraph correspondence principle]\label{hcp}  Fix $d \geq 2$.  For every $m \geq 1$, let $G^{(m)} = (V^{(m)}, E^{(m)})$ be a finite $d$-uniform hypergraph, and let $\P^{(m)}$ be as in Definition \ref{hypergraph-embed}.  Then there exists a subsequence $0 < m_1 < m_2 < \ldots$ of $m$, and a probability measure $\P^{(\infty)}$ on the hypergraph universal space $(\Omega,\B_{\max})$, such that we have the weak convergence property \eqref{weak-conv}
and the permutation invariance property \eqref{pshift-2}.
Furthermore, we have the following relative independence property: for any $I, I_1,\ldots, I_l \in \N$ with $I \cap I_1 \cap \ldots \cap I_l$ infinite, the factors $\B_I$ and $\bigvee_{i=1}^l \B_{I_i}$ are relatively independent conditioning on 
$\bigvee_{i=1}^l \B_{I \cap I_i}$, with respect to this probability measure $\P^{(\infty)}$.  
\end{theorem}

Similarly, by repeating the proof of Lemma \ref{irl} almost word for word we obtain

\begin{lemma}[Infinitary hypergraph regularity lemma]\label{hrl}  Fix $d \geq 2$, and let $\P$ be a probability measure on the hypergraph universal space $(\Omega, \B_\max)$ which is permutation-invariant in the sense of \eqref{pshift-2}.  Then for any $I, I_1,\ldots, I_l \in \N$ with $I \cap I_1 \cap \ldots \cap I_l$ infinite, the factors $\B_I$ and $\bigvee_{i=1}^l \B_{I_i}$ are relatively independent conditioning on  $\bigvee_{i=1}^l \B_{I \cap I_i}$.  
\end{lemma}

\begin{remark}
We should emphasise just how easily the regularity lemma has extended to the hypergraph case here.  This is contrast with the development of the finitary hypergraph regularity lemma, which has only been satisfactorily achieved quite recently \cite{nrs}, \cite{rs}, \cite{gowers-hyper}, \cite{tao-hyper} (with preliminary work in \cite{frankl}, \cite{chung-hyper}, \cite{frankl02}).  In the author's view this is because the regularity lemma is a relatively ``soft'' component of the theory; in the infinitary framework, the ``hard'' components of the theory are now isolated in 
the three fundamental extension
properties in Lemma \ref{weakmix}, Lemma \ref{finite}, Lemma \ref{zorn} (and to a lesser extent in Lemma \ref{repeat}). These three lemmas are roughly analogous to the ``counting lemma'' components of the hypergraph theory (although Lemma \ref{zorn} also captures some of the nature of the ``regularity lemma'' component, and is the step which is most responsible for the extremely poor quantitative bounds in this theory).  Unsurprisingly, it is also these three lemmas where one does the most non-trivial manipulation of small quantities such as $\eps$.  Fortunately, the infinitary setting allows one to isolate these epsilons from each other, despite the fact that all three of these basic lemmas are used repeatedly in the proof of the infinitary hypergraph removal lemma (Theorem \ref{compact-triangle}).  If instead we expanded out all of these lemmas within the proof of Theorem \ref{compact-triangle}, and allowed the various epsilons to mix together (with the order of quantifiers, etc. being carefully recorded), one would eventually end up with a complicated situation roughly analogous to those in the finitary proofs
\cite{nrs}, \cite{rs}, \cite{rodl}, \cite{rodl2}, \cite{gowers-hyper}, \cite{tao-hyper} of the hypergraph removal lemma.  Thus the infinitary perspective allows for a powerful \emph{encapsulation} of distinct components of the argument which greatly cleans up and clarifies the high-level structure of the proof, though the low-level components are, at a fundamental level, essentially the same as in the finitary approach.
\end{remark}

\section{An infinitary proof of the hypergraph removal lemma}

We can now repeat the arguments from Section \ref{triangle-sec} to obtain the following triangle-removal lemma of Nagle, Schacht, R\"odl, and Skokan
\cite{nrs}, \cite{rs}, \cite{rodl}, \cite{rodl2} (and independently by Gowers \cite{gowers-hyper}; see also \cite{tao-hyper} for a later proof):

\begin{theorem}[Hypergraph removal lemma]\label{hyper-remove} Fix $d \geq 2$, and let $G_0 = (V_0, E_0)$ be a $d$-uniform hypergraph.  Let $G = (V,E)$ be a $d$-uniform hypergraph with $|V|=n$ vertices.  Suppose that $G$ contains fewer than $\delta n^{|V_0|}$ copies of $G_0$ for some $0 < \delta \leq 1$, or more precisely
$$ | \{ (x_i)_{i \in V_0} \in V^{V_0}: \{ x_i: i \in e \} \in E \hbox{ for all } e \in E_0 \} | \leq \delta n^{|V_0|}.$$
Then it is possible to delete $o_{\delta \to 0;G_0,d}(n^d)$ edges from $G$ to create a $d$-uniform hypergraph $G'$ which has no copies of $G_0$ whatsoever.  Here the subscripting of the $o()$ notation by $G_0,d$ indicates that the quantity $o_{\delta \to 0;G_0}(n^d)$, when divided by $n^d$, goes to zero as $\delta \to 0$ for each fixed $G_0,d$, but the decay rate is not uniform in $G_0,d$.
\end{theorem}

\begin{remark}
As with the triangle removal lemma, this lemma has previously only been proven via a hypergraph regularity lemma, followed by a counting lemma.  This is rather complicated; the shortest proof known (in \cite{tao-hyper}) is about 25 pages, and requires some quite delicate computations.  While this current proof is arguably longer than the proof in \cite{tao-hyper}, and certainly less elementary, there are far fewer computations involved, and 
we believe the argument here is more conceptually clear.  This theorem has a number of applications, most notably in giving a proof not only of Szemer\'edi's theorem (Theorem \ref{sz-quant}) but also a multidimensional version due to Furstenberg and Katznelson \cite{fk}; see e.g. \cite{rodl2}, \cite{gowers-hyper}, \cite{tao-multiprime} for further discussion of this connection, and \cite{rstt} for some more applications of this theorem.  A variant of this theorem was also used in \cite{tao-multiprime} to establish that the Gaussian primes contain arbitrarily shaped constellations; we shall discuss this variant shortly.
\end{remark}

\begin{proof}(Sketch)  This is basically a repetition of the proof of Lemma \ref{trl}, so we sketch the main points only.  Fix $d$, $G_0$.
We can relabel $V_0$ to be $\{1,\ldots,n_0\}$ for some integer $n_0$; we can also easily assume that $E_0$ is non-empty.  If the theorem failed,
we argue much as in the proof of Lemma \ref{trl}, with $\{1,\ldots,n_0\}$ playing the role of $\{1,2,3\}$ (and thus $\{n_0+1,n_0+2,\ldots\}$
playing the role of $\{4,5,\ldots\}$).  We apply the hypergraph correspondence principle to obtain an infinitary limiting system
$(\Omega, \B_\max, \P^{(\infty)})$, and apply Theorem \ref{compact-triangle} with 
$J := \{1,\ldots,n_0\}$, $\ideal_\max := \{ e: e \subseteq e' \hbox{ for some } e' \in V_0 \}$, $E_e$ set equal to $A_e$ if
$e \in V_0$ and $E_e = \Omega$ otherwise (the latter happens precisely when $|e| < d$), and 
with $\B_e$ set equal to $\B_{e \cup \{n_0+1,n_0+2,\ldots\}}$ for all
$e \in \ideal_\max$.  One then continues the argument as in Lemma \ref{trl} (with the factor $100$ in \eqref{puff} replaced
by at least $2^d d! |E_0|$); the remainder of the proof proceeds with only the obvious minor changes.
\end{proof}

\begin{remark}  These results have analogues for partite hypergraphs (see \cite{tao-hyper}) and are proven similarly, but we will not do so here; the main difference is that instead of sampling all vertices from a single vertex class, one samples countably many vertices from each vertex class (which also leads to a more complicated symmetry group than $S_\infty$).
Just as the triangle removal lemma, Lemma \ref{trl}, has a stronger version in Lemma \ref{trs} which gives a complexity bound on
the approximating graph $G'$, the hypergraph removal lemma given above also comes with a stronger version, in which the approximating hypergraph
$G'$ is no longer a subhypergraph of $G$, but can be described using a partition of $\binom{V}{d-1}$ into $O_{\delta,G_0,d}(1)$ components.
We will neither state nor prove this stronger version here (the proof is much the same as Lemma \ref{trs}), but see \cite{tao-hyper} for an extremely similar statement (in the setting of partite hypergraphs rather than non-partite hypergraphs).  This version played an important role in the result in \cite{tao-multiprime} that the Gaussian primes contained arbitrarily shaped constellations.
\end{remark}
 
\appendix

\section{Review of probability theory}\label{probtools}

In this appendix we review the notation and tools from probability that we shall need.  There are two concepts here of particular importance: 
the concept of \emph{relative independence} of two or more factors in a probability space; and the ability to approximate complicated events or random variables by combinations of more elementary events or random variables.

\subsection{The algebra of events}

A probability space has two major structures; the set-theoretic structure of its events, and the measure-theoretic structure of the probability measure $\P$.  Because we will be dealing with multiple event spaces with a single probability measure, or multiple probability measures on a single event space, it will be conceptually clearer if we treat these two structures separately.  We begin with the structure of the event spaces.  For technical reasons it is convenient to restrict attention to countably generated spaces.

\begin{definition}[Event spaces]\label{event-def}  An \emph{event space} is a pair $(\Omega, \B_{\max})$, where the \emph{sample space} $\Omega$ is 
a non-empty set (possibly infinite), and $\B_{\max}$ is a $\sigma$-algebra on $\Omega$, i.e. a collection of subsets of $\Omega$ which are closed under countable unions, intersections, and complements, and which contains the empty set $\empty$ and $\Omega$.  We will also require that the $\sigma$-algebra $\B_{\max}$ be \emph{countably generated}, thus there exists a countable sequence of events $E_1, E_2, \ldots \in \B_\max$ such that $\B_\max$ is the minimal $\sigma$-algebra containing all these events.  We refer to elements of $\B_{\max}$ as \emph{(measurable) events}; we abuse notation and identify properties $P(x)$ of points $x \in \Omega$ with the associated event $\{ x \in \Omega: P(x) \hbox{ is true} \}$, and refer to the event simply as $P$.  If $A$ and $B$ are events, we use $A \vee B$ to denote the event that at least one of $A$ and $B$ are true (i.e. $A \vee B$ is the union of $A$ and $B$) and $A \wedge B$ to denote the event that $A$ and $B$ are not true (i.e. $A \wedge B$ is the intersection of $A$ and $B$).  We also use $\overline{A}$ to denote the event that $A$ is not true (thus $\overline{A} = \Omega \backslash A$).  
\end{definition}

\begin{example} If $\Omega$ is at most countable, the \emph{power-set event space} $(\Omega, 2^\Omega)$ of a set $\Omega$ is achieved by setting $\B_{\max} := 2^{\Omega} := \{ E: E \subseteq \Omega \}$ to be the power set of $\Omega$.  (If $\Omega$ is uncountable, $2^\Omega$ is no longer countably generated.)
\end{example}

\begin{definition}[Factors]\label{factor-def} Let $(\Omega,\B_{\max})$ be an event space.
A \emph{factor} is a subset $\B$ of $\B_{\max}$ which is also a countably generated $\sigma$-algebra.  More generally, we say that $\B_1$ is a \emph{factor of} $\B_2$ (or $\B_2$ \emph{extends} $\B_1$) if $\B_1, \B_2$ are both $\sigma$-algebras in $\B_{\max}$ and $\B_1 \subseteq \B_2$.  We say
that a factor is \emph{finite} if it consists of only finitely many events, thus for instance the \emph{trivial factor} $\{ \emptyset, \Omega\}$ is finite.   An event is \emph{$\B$-measurable} if it lies in $\B$.  A \emph{random variable} is any function $f: \Omega \to \R$ with the property that the events $f \in V$ are $\B_{\max}$-measurable for all open sets $V$; if these events are in fact $\B$-measurable, we say that the random variable $f$ is $\B$-measurable also.  In particular, if an event $E$ is \emph{$\B$-measurable}, then its \emph{indicator variable} $\I(E)$, defined to equal $1$ when $E$ is true and $0$ otherwise, is also $\B$-measurable.
If ${\mathcal E} \subseteq \B_{\max}$ is any collection of events, we let $\B[{\mathcal E}]$ denote the factor \emph{generated by} these events (i.e. the intersection of all factors that contain ${\mathcal E}$).  In particular, if $E$ is a single event, we let $\B[E] = \{ \emptyset, E, \overline{E}, \Omega \}$ denote the (finite) factor 
generated by $E$.  Similarly, if $X$ is a random variable taking finitely many values, we use $\B[X]$ to denote the factor generated by the events $X = c$, where $c$ ranges over the range of $X$.
We write $\B_1 \vee \B_2$ for $\B[\B_1 \cup \B_2]$, thus $\B_1 \vee \B_2$ is the least common extension of $\B_1$ and $\B_2$.  More generally, we can define the least common extension $\bigvee_{\alpha \in A} \B_\alpha$ of any at most countable collection of factors $\B_\alpha$.    
\end{definition}

\begin{example}[Finite factors] Let $A_1,\ldots,A_n$ be a partition of the sample space $\Omega$ into disjoint non-empty events.  Then $\B = \B[A_1,\ldots,A_n]$ is the finite factor consisting of all events which are the union of zero or more of the $A_i$ (and all finite factors are of this form).  We refer to $A_i$ as the \emph{atoms} of $\B$.  Let $i: \Omega \to \{1,\ldots,n\}$ be the random variable which indexes which atom one lies in, thus $x \in A_{i(x)}$ for all $x \in \Omega$. A random variable $f$ is $\B$-measurable if and only if it is determined by $i$, thus $f(x) = F(i(x))$ for some function $F: \{1,\ldots,n\} \to \R$.  One finite factor $\B_1$ extends another $\B_2$ if the partition into $\B_1$-atoms is finer than the partition into $\B_2$-atoms (thus every $\B_2$-atom is the union of $\B_1$-atoms).
\end{example}

We shall also need the notion of a \emph{(boolean) algebra}, namely a subset $\B$ of $\B_\max$ which is closed under \emph{finite} intersections, unions, complements, and contains $\emptyset$ and $\Omega$.  Thus every factor is an algebra, but not conversely.  The reason we need to deal with algebras rather than factors is because of the observation that the algebra generated by a countable sequence of events remains countable (indeed it is nothing more than the collection of finite boolean combinations of events from that sequence), whereas the factor generated by the same sequence can be uncountable.  This is important when applying the Arzela-Ascoli diagonalisation argument (see Lemma \ref{arzela} below).

\begin{example}  Let $\Omega = [0,1)^2$, and let $\B_{\max}$ be the Borel $\sigma$-algebra (i.e. the algebra generated by the open sets).  Let $\B_{\reg}$ be the space of \emph{elementary sets}, defined as the finite unions of half-open rectangles $[a,b) \times [c,d)$ where $a,b,c,d$ are rational.  Then $\B_{\reg}$ is an algebra but not a factor, and is countable; furthermore $\B_{\max}$ is generated by $\B_{\reg}$.
\end{example}

\subsection{Probability spaces}

We now add the structure of a probability measure to an event space, to form a probability space.

\begin{definition}[Probability spaces]  A \emph{probability space} is a triplet $(\Omega, \B_{\max}, \P)$, where $(\Omega,\B_{\max})$ is an event space, and $\P: \B_{\max} \to [0,1]$ is a \emph{probability measure}, i.e. a countably additive non-negative measure on  $\B_{\max}$ with $\P(\Omega) = 1$.  A \emph{null event} is an event of probability zero.  A statement is true \emph{almost surely} if it is only false on a null event.
\end{definition}

\begin{remark}
We do not assume our event space $(\Omega,\B_{\max})$ to be complete.  Thus, it is not necessarily the case that any subset of a null event is still a measurable event.  (It may help to think of the $\sigma$-algebras here as being like Borel $\sigma$-algebras --- that is, algebras generated by open sets --- rather than Lebesgue $\sigma$-algebras.)
\end{remark}

In the remainder of this appendix we assume that the probability space $(\Omega, \B_{\max}, \P)$ is fixed.

\begin{definition}[Random variables]  We consider two random variables equivalent if they are almost surely equal.  If $f$ is absolutely integrable, we use $\E(f)$ to denote the integral of $f$ with respect to the probability measure $\P$, and write $\|f\|_{L^1} = \|f\|_{L^1(\B_{\max};\P)}$ for $\E (|f|)$.  Thus for instance $\E(\I(E)) = \P(E)$ for any event $E$.  Similarly, we write $\|f\|_{L^2} = \|f\|_{L^2(\B_{\max};\P)}$ for $\E(|f|^2)^{1/2}$ whenever $f$ is square-integrable, and $\|f\|_{L^\infty} = \|f\|_{L^\infty(\B_{\max};\P)}$ for the essential supremum of $f$.  We will drop the measure $\P$, and sometimes the factor $\B_{\max}$, from the $L^p(\B_{\max};\P)$ notation when these are clear from context.
\end{definition}

It will be important to develop relative versions of all these concepts with respect to factors of $\B_{\max}$.

\begin{definition}[Conditional expectation] If $p=1,2,\infty$ and $\B$ is a factor, we let $L^p(\B) = L^p(\B;\P)$ denote the space of $\B$-measurable random variables with finite $L^p$ norm (identifying two random variables if they are equivalent).  Observe that $L^2(\B)$ is a Hilbert space with inner product $\langle f, g \rangle := \E(fg)$; since $\B$ is countably generated, we see that $L^2(\B)$ is separable. We define the \emph{conditional expectation operator}
$f \mapsto \E(f|\B)$ to be the orthogonal projection from $L^2(\B_{\max})$ to $L^2(\B)$; note that $\E(f|\B)$ is only defined up to almost sure equivalence.  If $E$ is an event, we write $\P(E|\B)$ for $\E(\I(E)|\B)$, and refer to $\P(E|\B)$ as the \emph{conditional probability} of $E$ with respect to the factor $\B$.    
\end{definition}

We have the useful

\begin{lemma}[Pythagoras' theorem]\label{pythag}  Let $\B'$ be an extension of $\B$.  Then for any $f \in L^2(\B_\max)$ we have
$$ \| \E(f|\B') \|_{L^2}^2 = \| \E(f|\B) \|_{L^2}^2 + \| \E(f|\B') - \E(f|\B) \|_{L^2}^2.$$
\end{lemma}

\begin{proof}  This follows since $\E(f|\B')$ is the orthogonal projection to $L^2(\B')$, and $\E(f|\B)$ is the orthogonal projection to the smaller space $L^2(\B)$.
\end{proof}

\begin{remark} In this paper we shall deal almost exclusively with bounded random variables (indeed, they will almost always be bounded between $-1$ and $1$).  Thus issues of \emph{integrability} will not be a concern to us; this also means that we do not have to distinguish between convergence in $L^1$, convergence in $L^2$, and convergence in measure.  It will however be crucial to keep track the \emph{measurability} of our random variables with respect to the various factors involved in the argument.
\end{remark}

\begin{example}[Finite factors] Let $\B$ be a finite factor with atoms $A_1,\ldots,A_n$.  If $f \in L^\infty(\B_{\max})$, the conditional expectation $\E(f|\B)$ is well-defined on all atoms $A_i$ of non-zero probability, and is equal to $\E(f|A_i) := \E( f\I(A_i) ) / \P(A_i) )$ on each such atom.  Similarly we have $\P(E|\B) = \P( E|A_i) := \P(E \wedge A_i ) / \P(A_i)$ on such atoms.  Of course one can develop similar explicit formulae for the conditional covariance of two random variables or events.
\end{example}

We recall some very standard properties of conditional expectation, that we shall use without further comment.  
The conditional expectation operation $f \mapsto \E(f|\B)$ is linear, positivity preserving, and is
a contraction on $L^p$ for $p=1,2,\infty$.  In particular conditional expectation is continuous in each of the $L^p$ topologies, which allows us to easily apply density arguments when verifying identities involving conditional expectation (i.e. it suffices to verify such identities for a dense subclass of random variables, such as simple random variables).
We also have the module property that $\E(fg|\B) = f \E(g|\B)$ whenever $f \in L^\infty(\B)$ and $g \in L^\infty(\B_{\max})$.  

In order to pass from a sequence of finitary objects to an infinitary one, the following lemma will be crucial.

\begin{lemma}[Arzela-Ascoli diagonalisation argument]\label{arzela}  Let $\P^{(1)}, \P^{(2)}, \ldots$ be a sequence of probability measures on an event space $(\Omega, \B_{\max})$. Let $\B_{\reg}$ be a countable algebra which generates $\B_{\max}$ as a $\sigma$-algebra.  Then there exists a subsequence $0 < k_1 < k_2 < \ldots$ of integers and a probability measure $\P$ such that 
$$ \lim_{i \to \infty} \P^{(k_i)}(F) = \P(F) \hbox{ for all } F \in \B_\reg.$$
In other words, $\P^{(k_i)}$ is weakly convergent to $\P$, when tested against the algebra of events $\B_\reg$.
\end{lemma}

\begin{proof} We enumerate $\B_\reg$ as $F_1, F_2, \ldots$ (duplicating events if necessary, if $\B_\reg$ happens to be finite).  By using the sequential compactness of the unit interval $[0,1]$ (i.e. the Heine-Borel theorem), we can obtain a 
sequence $k_{1,1} < k_{1,2} < k_{1,3} < \ldots$ such that $\P^{(k_{1,i})}(F_1)$ converges as $i \to \infty$ to a limit, say $p_1 \in [0,1]$.  Then we can extract a subsequence $k_{2,1} < k_{2,2} < k_{2,3} < \ldots$ of that sequence such that $\P^{(k_{2,i})}(F_2)$ converges as $i \to \infty$ to a limit, say $p_2 \in [0,1]$.  We continue in this fashion and then extract the diagonal sequence $k_i := k_{i,i}$ to obtain a sequence $p_1,p_2,\ldots \in [0,1]$ such that $\lim_{i \to \infty} \P^{(k_i)}(F_j) = p_j$ for each $j =1,2,\ldots$.  One easily verifies that the map $F_j \mapsto p_j$ is \emph{finitely} additive, non-negative, and maps $\emptyset$ to $0$ and $\Omega$ to $1$.  Invoking the Kolmogorov extension theorem (or the Carath\'eodory extension theorem) we can construct a probability measure $\P$ such that $\P(F_j) = p_j$, and the claim follows.
\end{proof}

\begin{remark}\label{construct} One can also obtain this lemma from the Banach-Alaoglu theorem and the Riesz representation theorem (though one should take care to distinguish the notions of compactness and sequential compactness).  Observe that both the Heine-Borel theorem and the Kolmogorov extension theorem are completely constructive, and so this lemma does not use the axiom of choice.  See \cite{tao:transference} for further discussion.
\end{remark}

\subsection{Approximation lemmas}

We will frequently need to approximate a random variable or event in a complicated factor by linear, polynomial, or boolean combinations of random variables or events in simpler factors.  To do this we shall use some very simple and standard tools, which we collect here for the reader's convenience.

Recall that a random variable is \emph{simple} if it only takes on finitely many values, or equivalently if it is the finite linear combination of
indicator functions, or equivalently if it is measurable with respect to a finite factor.  
The following lemma is standard in measure theory:

\begin{lemma}\label{dense} Let $\B$ be a factor and $p=1,2,\infty$.  Then the simple random variables in $L^p(\B)$ are dense in $L^p(\B)$.  
\end{lemma}

Because of this, the task of approximating random variables quickly boils down to approximating events.
Let us say that two events $E, F$ are \emph{$\eps$-close} if $\P( E \backslash F ) + \P( F \backslash E) \leq \eps$.  

\begin{lemma}[Approximation by finite complexity events]\label{product-approx} Let $\B = \B[{\mathcal E}]$ be a factor generated by a (possibly infinite) collection $\mathcal{E}$ of events, and let $\eps > 0$.  Then every event in $\B$ is $\eps$-close to a finite boolean combination of events from $\mathcal{E}$.    In particular, if $\B$ is generated by an algebra $\B_{\reg}$, then every event in $\B$ is $\eps$-close to an event from $\B_\reg$.  If $f \in L^1(\B)$, then there exists a finite factor $\B'$ of $\B$ generated by finitely many events in $\mathcal{E}$, such that
$\|f - \E(f|\B')\|_{L^1(\B)} \leq \eps$.
\end{lemma}

\begin{proof}  Let $\B_{\reg}$ be the algebra generated by $\mathcal{E}$ (i.e. the space of finite boolean combinations of events from $\mathcal{E}$). Let $\B_\eps$ denote the collection of events which is $\eps$-close to an element of $\B_{\reg}$.  Then one easily verifies that $\bigcap_{\eps > 0} \B_\eps$ is a factor that contains $\mathcal{E}$, and thus contains $\B$, and the first and second claims follow.  To prove the final claim, first use Lemma \ref{dense} to reduce to the case where $f$ is simple, and then use linearity to reduce to the case when $f=\I(E)$ is an indicator function.  By the previous claims, we can find an event $E' \in \B_{\reg}$ which is $\eps/2$-close to $E$, thus $\| f - \I(E') \|_{L^1(\B)} \leq \eps/2$.  This $E'$ lies in some finite factor $\B'$ generated by $\mathcal{E}$, and thus on taking conditional expectations in $\B'$ we see that $\| \E(f|\B') - \I(E') \|_{L^1(\B)} \leq \eps/2$.  The claim now follows from the triangle inequality.
\end{proof}

\begin{corollary}[Limits of chains]\label{product-monotone}  Let $A$ be a totally ordered set, and for each $\alpha \in A$ let $\B_\alpha$ be a factor of $\B_{\max}$ with the monotonicity property $\B_\alpha \subseteq \B_\beta$ whenever $\alpha \leq \beta$.  Let $\B := \bigvee_{\alpha \in A} \B_\alpha$.  Then for any $f \in L^2(\B)$, the net $\E(f|\B_\alpha)$ converges to $f$ in $L^2$ norm (thus for every $\eps > 0$ there exists $\beta \in A$ such that
$\|f - \E(f|\B_\alpha)\|_{L^2(\B)} \leq \eps$ for all $\alpha \geq \beta$).
\end{corollary}

\begin{proof} Let $\eps > 0$.  Applying Lemma \ref{product-approx} with ${\mathcal E} = \bigcup_{\alpha \in A} \B_\alpha$, we can find a finite factor $\B'$ generated by finitely many events in $\mathcal{E}$ such that $\|f - \E(f|\B') \|_{L^2(\B)} \leq \eps$. By monotonicity we see that $\B'$ is a factor of $\B_\alpha$ for some $\alpha \in A$. The claim then follows from Pythagoras' theorem.
\end{proof}

\begin{corollary}[Approximation by finite factors]\label{product-2} Let $\B_1,\ldots,\B_k$ be factors and $\eps > 0$.  Then every event in $\B_1 \vee \ldots \vee \B_k$ is $\eps$-close to a finite boolean combination of events in $\B_1 \cup \ldots \cup \B_k$.  Furthermore, given any random variable $f \in L^\infty(\B_1 \vee \ldots \vee \B_k)$, there exists finite factors $\B'_i$ of $\B_i$ for $i=1,\ldots,k$ respectively such that 
$\| f- \E(f|\B'_1 \vee \ldots \B'_k) \|_{L^1(\B_1 \vee \ldots \vee \B_k)} \leq \eps$.
\end{corollary}

\begin{proof} The first claim follows from Lemma \ref{product-approx} by setting ${\mathcal E} := \B_1 \cup \ldots \cup \B_k$.  To verify the second claim, first use Lemma \ref{product-approx} to locate a finite factor $\B'$ generated by finitely many elements in $\B_1 \cup \ldots \cup \B_k$ such that
$\|f - \E(f|\B') \|_{L^1(\B_1 \vee \ldots \vee \B_k)} \leq \eps/2$.  Now observe that $\B'$ is a factor of $\B'_1 \vee \ldots \vee \B'_k$ for some finite factors $\B'_i$ of $\B_i$ for $i=1,\ldots,k$.  The claim now follows from the same triangle inequality argument used to prove Lemma \ref{product-approx}.
\end{proof}

\subsection{Relative independence}

Now we come to a fundamental notion for us, namely that of (relative) independence of two or more factors.

\begin{definition}[Independence]  We say that two factors $\B_1, \B_2$ are \emph{unconditionally independent} if we have 
$$\E(f_1 f_2) = \E(f_1) \E(f_2)$$
for all $f_1 \in L^\infty(\B_1)$ and $f_2 \in L^\infty(\B_2)$.  More generally, we say that two factors $\B_1, \B_2$ are \emph{relatively independent} conditioning on a third factor $\B$ with respect to the probability measure $\P$ if we have 
\begin{equation}\label{eff-0}
 \E( f_1 f_2 | \B ) = \E(f_1|\B) \E(f_2|\B)
 \end{equation}
almost surely for all $f_1 \in L^\infty(\B_1)$ and $f_2 \in L^\infty(\B_2)$.  In many cases, the probability measure $\P$ will be clear from context and we shall omit the phrase ``with respect to $\P$''. Given an at most countable collection of factors $(\B_\alpha)_{\alpha \in A}$, we say that these factors are \emph{jointly unconditionally independent} (resp. \emph{jointly relatively independent} conditioning on a factor $\B$) if $\bigvee_{\alpha \in A_1} \B_\alpha$ and $\bigvee_{\alpha \in A_2} \B_\alpha$ are unconditionally independent (resp. relatively independent conditioning on $\B$) for all disjoint subsets $A_1,A_2$ of $A$.  We say that a collection of events $E_1, E_2, \ldots$ is unconditionally independent (resp. relatively independent conditioning on $\B$) if their associated factors $\B[E_1], \B[E_2], \ldots$ are unconditionally independent (resp. relatively independent conditioning on $\B$).
\end{definition}

\begin{examples} Two factors $\B_1,\B_2$ are unconditionally independent if and only if $\P(E \wedge F) = \P(E) \P(F)$ for all $E \in \B_1$ and $F \in \B_2$.  In particular, two events $E$ and $F$ are unconditionally independent if and only if $\P(E \wedge F) = \P( E ) \P(F)$.
Three factors $\B_1,\B_2,\B_3$ are jointly unconditionally independent if and only if $\P( E_1 \wedge E_2 \wedge E_3 ) = \P(E_1) \P(E_2) \P(E_3)$ for all $E \in \B_1$ and $F \in \B_2$.  On the other hand, in order for three events $E,F,G$ to be jointly independent it is not quite enough that
$\P(E \wedge F \wedge G) = \P(E) \P(F) \P(G)$; one also needs $E,F,G$ to be pairwise independent, thus for instance $\P(E \wedge F) = \P(E) \P(F)$.
If $\B_1,\B_2,\B_3$ are jointly unconditionally independent, then $\B_1 \vee \B_3$ and $\B_2 \vee \B_3$ are conditionally independent over $\B_3$, even though they are almost certainly not unconditionally independent.  On the other hand, $\B_1$ and $\B_2$ are both unconditionally independent, and conditionally independent over $\B_3$.
\end{examples}

\begin{example}  Let $x_1,x_2,x_3$ be three elements chosen uniformly and independently at random from $\{0,1\}$.  Then the events $x_1=x_3$ and $x_2=x_3$ are unconditionally independent, but they are not relatively independent conditioning on the factor $\B[x_1=x_2]$.  Thus we see that unconditional independence is neither stronger nor weaker than relative independence.  
\end{example}

Taking expectations in \eqref{eff-0} we obtain
\begin{equation}\label{eff}
\E( f_1 f_2 ) = \E( \E(f_1|\B) f_2 ) = \E( f_1 \E(f_2|\B) ) = \E( \E(f_1|\B) \E(f_2|\B) )
\end{equation}
whenever $\B_1,\B_2$ are relatively independent conditioning on $\B$, $f_1 \in L^\infty(\B_1)$ and $f_2 \in L^\infty(\B_2)$.

There are several equivalent formulations of relative independence.

\begin{lemma}\label{verify-independence}  Let $\B_1, \B_2, \B$ be factors.  Then the following are equivalent:
\begin{itemize}
\item[(i)] $\B_1$ and $\B_2$ are relatively independent conditioning on $\B$.
\item[(ii)] We have $\E(f_1 | \B \vee \B_2) = \E(f_1|\B)$ almost surely for all $f_1 \in L^2(\B_1)$.
\item[(iii)] We have $\| \E( f_1 | \B \vee \B_2 )\|_{L^2} = \|\E( f_1 | \B )\|_{L^2}$ for all $f_1 \in L^2(\B_1)$.
\item[(iv)] We have $\| \P( E_1 | \B \vee \B_2 )\|_{L^2} = \|\P( E_1 | \B )\|_{L^2}$ for all $E_1 \in \B_1$.
\end{itemize}
\end{lemma}

\begin{proof}  The equivalence of (ii) and (iii) follows from Lemma \ref{pythag}.  The equivalence of (iii) and (iv)
follows from Lemma \ref{dense}, linearity, and a standard limiting arguments.  

To see that (ii) implies (i), observe for $f_1 \in L^\infty(\B_1)$ and $f_2 \in L^\infty(\B_2)$ that
\begin{align*}
\E( f_1 f_2 | \B ) &= \E( \E(f_1 f_2|\B \vee \B_2) | \B ) \\
&=\E( \E(f_1|\B \vee \B_2) f_2 | \B ) \\
&=\E( \E(f_1|\B) f_2 | \B ) \\
&=\E(f_1|\B) \E( f_2 | \B ) 
\end{align*}
where we have used the module property twice.

Finally, we show that (i) implies (iv).  We observe from (i) and the module property that
$$ \E( f_1 f_2 h | \B ) = \E( f_1|\B) \E( f_2 h | \B )$$
whenever $f_1 \in L^\infty(\B_1)$, $f_2 \in L^\infty(\B_2)$, and $h \in L^\infty(\B)$.  Taking linear combinations 
and using limiting arguments we conclude that
$$ \E( f_1 g | \B ) = \E( f_1|\B) \E( g | \B )$$
whenever $g \in L^\infty(\B \vee \B_2)$.  We take expectations and obtain
$$ \E(f_1 g) = \E(\E(f_1|\B) \E(g|\B)).$$
Applying this with $f_1 := \I(E_1)$ and $g = \P(E_1| \B \vee \B_2)$ we obtain
$$ \|\P( E_1 | \B \vee \B_2 )\|_{L^2}^2 = \E(\I(E_1) \P( E_1 | \B \vee \B_2 ))
= \E( \P(E_1|\B) \E( \P( E_1 | \B \vee \B_2 ) | \B ) ) = \| \P(E_1|\B) \|_{L^2}^2$$
and (iv) follows.
\end{proof}

Now we can observe the following stability properties concerning relative independence.

\begin{proposition}\label{stability}  Let $\B_1, \B_2$ be two factors which are relatively independent conditioning on another factor $\B$.
\begin{itemize}
\item[(i)] (Monotonicity) If $\B'_1$ is a factor of $\B_1$ and $\B'_2$ is a factor of $\B_2$, then $\B'_1$ and $\B'_2$ are relatively independent  conditioning on $\B$.
\item[(ii)] (Absorption) $\B_1 \vee \B$ and $\B_2 \vee \B$ are relatively independent conditioning on $\B$.
\item[(iii)] (Gluing) Let $\B_3$ be a $\sigma$-algebra which is relatively independent of $\B_1 \vee \B_2$ conditioning on $\B$.  Then $\B_1$ is relatively independent of $\B_2 \vee \B_3$ conditioning on $\B$.
\item[(iv)] (Factors do not affect relative independence) If $\B'_1$ is a factor of $\B_1$ and $\B'_2$ is a factor of $\B_2$, 
then $\B_1$ and $\B_2$ are relatively independent conditioning on $\B \vee \B'_1 \vee \B'_2$.
\item[(v)] (Independent information does not affect relative independence) Let $\B_3$ be a $\sigma$-algebra which is independent of $\B \vee \B_1 \vee \B_2$.  Then $\B_1$ is relatively independent of $\B_2 \vee \B_3$ conditioning on $\B$.
\end{itemize}
\end{proposition}

\begin{proof}  The claim (i) is trivial.  To prove (ii), observe from symmetry and iteration that it suffices to show that $\B_1 \vee \B$ and $\B_2$ are relatively independent conditioning on $\B$.  But this follows from two applications of Lemma \ref{verify-independence}.

To prove (iii), it suffices by Lemma \ref{verify-independence} (and symmetry) to show that
$$ \E( h | \B \vee \B_1 ) = \E( h | \B )$$
for all $h \in L^\infty(\B_2 \vee \B_3)$.  By density it suffices to show that
$$ \E( f_2 f_3 | \B \vee \B_1 ) = \E( f_2 f_3 | \B )$$
for all $f_2 \in L^\infty(\B_2)$ and $f_3 \in L^\infty(\B_3)$.  But this follows from the relative independence
hypotheses and the module property:
\begin{align*}
\E( f_2 f_3 | \B \vee \B_1 ) &= \E( \E(f_2 f_3|\B \vee \B_1 \vee \B_2) | \B \vee \B_1 ) \\
&= \E( f_2 \E( f_3 | \B \vee \B_1 \vee \B_2 ) | \B \vee \B_1 ) \\
&= \E( f_2 \E( f_3 | \B ) | \B \vee \B_1 ) \\
&= \E( f_2 | \B \vee \B_1 ) \E(f_3 | \B ) \\
&= \E( f_2 | \B ) \E( f_3 | \B ) \\
&= \E( f_2 f_3 | \B ).
\end{align*}

Now we prove (iv).  By symmetry and iteration it will suffice to show that
$\B_1$ and $\B_2$ are relatively independent conditioning on $\B \vee \B'_1$.  From Lemma \ref{verify-independence}
we already have
$$ \| E(f_2|\B) \|_{L^2} = \| \E(f_2|\B \vee \B_1)\|_{L^2}$$
for all $f_2 \in L^2(\B_2)$.  From Lemma \ref{pythag} we conclude
$$ \| E(f_2|\B) \|_{L^2} = \|\E(f_2|\B \vee \B'_1)\|_{L^2} = \| \E(f_2|\B \vee \B_1)\|_{L^2}$$
and the claim follows from another application of Lemma \ref{verify-independence}.

Finally, we prove (v).  If $\B_3$ is independent of $\B \vee \B_1 \vee \B_2$, then by the monotonicity and factor properties (i), (iv) we conclude that $\B_3$ is relatively independent of $\B_1 \vee \B_2$ conditioning on $\B$.  The claim (v) then follows from the gluing property (iii).
\end{proof}

\section{Connection with recurrence theorems}\label{conrecur}

We have just seen how infinitary probabilistic statements such as Theorem \ref{compact-triangle} can imply finitary graph statements such as
Lemma \ref{trl}; later we shall see that one can also deduce finitary hypergraph statements in this manner.  It is also well known
(see \cite{rsz}, \cite{frankl}, \cite{frankl02}, \cite{rodl}, \cite{rodl2}, \cite{gowers-hyper}, \cite{soly-roth}, \cite{tao-multiprime})
that these graph and hypergraph statements can in turn be used to deduce density results such as Szemer\'edi's theorem.  This in turn is known
by the Furstenberg correspondence principle to be equivalent to results such as the Furstenberg recurrence theorem.  Concatenating all these implications, one thus expects results such as Theorem \ref{compact-triangle} to be capable of implying results such as Theorem \ref{furst-thm} directly, without the need to pass back and forth between the infinitary and finitary settings.

Somewhat surprisingly, it appears to be somewhat difficult to achieve this goal; the best the author was able to do was simply to compose the
various implications discussed above to obtain a connection.  For sake of completeness, we sketch a special case of this connection here, but it is puzzling that there seems to be little ``synergy'' between these two infinitary results, despite their similarity.  As there appear to be no major new features emerging in this connection, we will skip over some of the details.

One can demonstrate the connection using the Furstenberg recurrence theorem (Theorem \ref{furst-thm}), but it will be slightly more convenient to instead work with the following variant:

\begin{theorem}[Furstenberg-Katznelson recurrence theorem, special case]\label{fk-thm}\cite{fk}  Let $(\Omega, \B_{\max}, \P)$ be a probability space.  Let $S,T:\Omega \to \Omega$ be two commuting probability-preserving bi-measurable maps.  Then for all events $A \in \B_{\max}$ with $\P(A) > 0$, we have
$$ \liminf_{N \to \infty} \frac{1}{2N+1} \sum_{n=-N}^N \P( A \wedge T^n A \wedge S^n A ) > 0.$$
\end{theorem}

This theorem is equivalent to the assertion that any subset of $\Z^2$ with positive upper density contains infinitely many right-angled
triangles $(x,y), (x+r,y), (x,y+r)$, a result first obtained by Ajtai and Szemer\'edi \cite{AS}.  In \cite{soly-roth} it was observed that this theorem
followed from the triangle removal lemma.  Setting $S := T^2$ we obtain the $k=3$ case of Theorem \ref{furst-thm}.  The full version of the Furstenberg-Katznelson recurrence theorem allows for an arbitrary number of commuting shifts, and can be treated by a modification of the arguments
presented here.

To transfer this theorem to a setting where Theorem \ref{compact-triangle} is applicable, we will have to perform essentially the entire machinery used in the graph correspondence principle.  It is convenient not to work with graphs on $\N$, but rather on tripartite graphs connecting three copies of $\Z$:

\begin{definition}[Tripartite graph universal space]\label{trigraph-universal}  A \emph{tripartite infinite graph} is a sextuple
$G = (\Z,\Z,\Z,E_{12},E_{23},E_{31})$ where $E_{12},E_{23},E_{31}$ are subsets of $\Z^2$. Let $\Omega^\Delta$ denote the space of all tripartite infinite graphs.  On this space $\Omega^\Delta$, we introduce the events $A_{ij,k_i,k_j}$ for $ij=12,23,31$ and $k_i,k_j \in \Z$ by $A_{ij,k_i,k_j} := \{ G \in \Omega: (k_i,k_j) \in E_{ij} \}$, and let $\B^\Delta_\max$ be the $\sigma$-algebra generated by the $A_{ij,k_i,k_j}$.  We also introduce the regular algebra $\B^\Delta_\reg$ generated by the $A_{ij,k_i,k_j}$.  For any three permutations $\sigma_1,\sigma_2,\sigma_3:\Z \to \Z$ we can define an action of $(\sigma_1,\sigma_2,\sigma_3)$ on $\B^\Delta_\max$ by mapping $A_{ij,k_i,k_j}$ to $A_{ij,\sigma_i(k_i),\sigma_j(k_j)}$.  For any subsets $I_1,I_2,I_3$
of $\N$, we define $\B^\Delta_I$ to be the factor of $\B^\Delta_\max$ generated by the events $A_{ij,k_i,k_j}$ where $ij=12,23,31$, $k_i \in I_i$, and $k_j \in I_j$.
\end{definition}

Now we embed the system in Theorem \ref{fk-thm} into this universal space.

\begin{definition}[Tripartite graph universal embedding]\label{trigraph-embed}  Let $(\Omega, \B_\max, \P)$, $S$, $T$, $A$ be as in Theorem \ref{fk-thm}.
Let $N \geq 1$ be a natural number. We introduce the probability space $(\Omega^{(N)}, \B_\max^{(N)}, \P^{(N)})$, defined as the space associated to sampling three infinite sequences $(n_{i,k_i})_{k_i \in \Z}$ for $i=1,2,3$ uniformly and independently at random from $[N] = \{1,\ldots,N\}$.  Thus
the product space $(\Omega \times \Omega^{(N)}, \B_\max \times \B_\max^{(N)}, \P \times \P^{(N)})$ represents the independent sampling of
a point $x$ from $\Omega$, together the three sequences $n_{i,k_i} \in [N]$ for $i=1,2,3$ and $k_i \in \Z$.  For any such $x$ and $n_{i,k}$,
we associate an infinite tripartite graph $G = (\Z,\Z,\Z, E_{12}, E_{23}, E_{31})$ in $\Omega^\Delta$ by setting
\begin{align*}
E_{12} &:= \{ (k_1,k_2) \in \Z \times \Z: T^{n_{1,k_1}} S^{n_{2,k_2}} x \in A \} \\
E_{23} &:= \{ (k_2,k_3) \in \Z \times \Z: T^{n_{3,k_3}-n_{1,k_1}} S^{n_{2,k_2}} x \in A \} \\
E_{31} &:= \{ (k_3,k_1) \in \Z \times \Z: T^{n_{1,k_1}} S^{n_{3,k_3}-n_{2,k_2}} x \in A \}.
\end{align*}
This is a measurable map from $\Omega \times \Omega^{(N)}$ to $\Omega^\Delta$ (the inverse image of $A_{12,k_1,k_2}$ is the measurable
event $T^{n_{1,k_1}} S^{n_{2,k_2}} x \in A$, and similarly for the other two classes of generating events), and so we can push forward
the measure $\P\times \P^{(N)}$ to a measure on $(\Omega \times \Omega^{(N)} \times \Omega^\Delta, \B_\max \times \B_\max^{(N)} \times \B_\max^\Delta)$, which by abuse of notation we will also call $\P\times \P^{(N)}$.
\end{definition}

A computation (using the probability-preserving and commuting nature of $T$ and $S$) shows that
\begin{align*}
 \P \times \P^{(N)}( &A_{12,0,0} \wedge A_{23,0,0} \wedge A_{31,0,0} )\\
&= \frac{1}{N^3} \sum_{n_1,n_2,n_3 \in [N]} \P( T^{n_1} S^{n_2} A \wedge T^{n_3-n_2} S^{n_2} A \wedge T^{n_1} S^{n_3-n_2} A ) \\
&= \frac{1}{N^3} \sum_{n_1,n_2,n_3 \in [N]} \P( A \wedge T^{n_3-n_2-n_1} A \wedge S^{n_3-n_2-n_1} A ) \\
&\leq \frac{1}{N} \sum_{n=-2N}^N \P( A \wedge T^n A \wedge S^n A ).
\end{align*}
Thus to prove Theorem \ref{fk-thm} it will suffice to show that
$$ \liminf_{N \to \infty}  \P \times \P^{(N)}( A_{12,0,0} \wedge A_{23,0,0} \wedge A_{31,0,0} ) > 0.$$
Suppose this were false.  Then one can find a sequence $N^{(m)}$ of $N$, going to infinity as $m \to \infty$, such that
$$ \lim_{m \to \infty} \P \times \P^{(N^{(m)})}( A_{12,0,0} \wedge A_{23,0,0} \wedge A_{31,0,0} ) = 0.$$
By applying Lemma \ref{arzela}, we can pass to a subsequence if necessary, and obtain a limiting probability measure
$\P^\Delta$ on $(\Omega^\Delta, \B_\max^\Delta)$ with the weak convergence property
$$ \lim_{m \to \infty} \P \times \P^{(N^{(m)})}( E ) = \P^\Delta(E) \hbox{ for all } E \in \B_\reg^\Delta.$$
The individual measures $\P \times\P^{(N)}$ can be easily verified to be invariant under triple permutations $(\sigma_1,\sigma_2,\sigma_3)$,
and so the limiting measure $\P^\Delta$ is also.

By adapting the arguments used to prove Lemma \ref{irl}, one can exploit the above invariance to 
establish the following relative independence property:
If $I_i, I_{i,1},\ldots, I_{i,l}$ are subsets of $\Z$ with $I_i \cap I_{i,1} \cap \ldots \cap I_{i,l}$ for $i=1,2,3$ , then the factors $\B_{I_1,I_2,I_3}$ and $\bigvee_{j=1}^l \B_{I_{1,j}, I_{2,j}, I_{3,j}}$ are relatively independent conditioning on $\bigvee_{j=1}^l \B_{I_1 \cap I_{1,j}, I_2 \cap I_{2,j}, I_3 \cap I_{3,j}}$, with 
respect to this probability measure $\P^\Delta$.  We omit the details of this as they are essentially the same as in the proof of Lemma \ref{irl} except for minor notational complications.

We now apply Theorem \ref{compact-triangle} on $(\Omega^\Delta, \B_\max^\Delta, \P^\Delta)$, 
with $J := \{1,2,3\}$, $\ideal_\max := \{ e: |e| \leq 2 \}$, $E_e$ set equal to $A_{ij,0,0}$ if
$e = \{i,j\}$ for some $ij=12,23,31$, and $E_e = \Omega^\Delta$ otherwise, and with $\B_e$ set equal to $\B_{I_1,I_2,I_3}$, where
$I_i$ is equal to $\Z$ if $i \in e$ and $\Z \backslash \{0\}$ if $i \not \in e$.  Thus for instance $\B_{\{1,2\}} = \B_{\Z,\Z,\Z \backslash \{0\}}$.
The hypotheses of the theorem are easily verified, and by arguing as in the proof of Lemma \ref{trl} (or Lemma \ref{trs}) we can find regular events
$F'_{i,j} \in \B_{\{i,j\}}$ for $ij=12,23,31$ obeying \eqref{fff2} and
$$   \P^{\Delta}( A_{i,j} \backslash F'_{i,j} ) < \eta / 10 \hbox{ for } ij = 12, 23, 31,$$
where $\eta := \P(A) > 0$ is the probability of the original event $A$.  In particular, for all sufficiently large $m$ we have
$$ \P\times\P^{N^{(m)}}( A_{i,j} \backslash F'_{i,j} ) < \eta / 10 \hbox{ for } ij = 12, 23, 31.$$
Now recall that the variables $n_{1,0}$, $n_{2,0}$, $n_{3,0}$ are independently and uniformly distributed on the interval $[N^{(m)}]$.
For any fixed $n_{1,0}, n_{2,0}$, the probability that $n_{3,0}$ equals $n_{1,0}+n_{2,0}$ will equal $1/N^{(m)}$ approximately half the time, and
0 the other half of the time.  Since the event $A_{1,2} \backslash F'_{1,2}$ is independent of $n_{3,0}$, we thus conclude from Bayes' formula that
$$ \P\times\P^{N^{(m)}}( A_{1,2} \backslash F'_{1,2} | n_{3,0} = n_{1,0} + n_{2,0} ) < \eta / 5$$
for $m$ sufficiently large.  Similar arguments in fact give
$$ \P\times\P^{N^{(m)}}( A_{i,j} \backslash F'_{i,j} | n_{3,0} = n_{1,0} + n_{2,0} ) < \eta / 5 \hbox{ for } ij=12,23,31.$$
On the other hand, from \eqref{fff2} we have\footnote{Note how important it is here that the event $F'_{1,2} \wedge F'_{2,3} \wedge F'_{3,1}$ be \emph{empty}, rather than merely being a null event with respect to $\P^\Delta$.  In the latter case, the event would have a small but nonzero measure in $\P \times \P^{N^{(m)}}$, and we would be unable to condition this event to the vanishingly small probability event $n_{3,0} = n_{1,0}+n_{2,0}$ without losing control on the conditional probability.  The point is that the constraint $n_{3,0} = n_{1,0}+n_{2,0}$ creates a ``diagonal measure'' which is singular with respect to $\P^\Delta$, and so null events in $\P^\Delta$ do not necessarily restrict to null events on the diagonal measure. However, events which have empty intersection with respect to $\P^\Delta$ will clearly continue to have empty intersection with respect to the diagonal measure.  This robustness with respect to change of measure is what makes Theorem \ref{compact-triangle} (which is basically a mechanism for converting null events to empty events) so powerful.} 
$$ \P\times\P^{N^{(m)}}( F'_{1,2}\wedge F'_{2,3}\wedge F'_{3,1} | n_{3,0} = n_{1,0} + n_{2,0} ) = 0.$$
Combining this with the preceding estimate we see that
$$ \P\times\P^{N^{(m)}}( A_{1,2}\wedge A_{2,3}\wedge A_{3,1} | n_{3,0} = n_{1,0} + n_{2,0} ) < 3 \eta/5.$$
However, the left-hand side equals
$$\P^{N^{(m)}}(  \P( T^{n_{1,0}} S^{n_{2,0}} A \wedge T^{n_{3,0}-n_{2,0}} S^{n_{2,0}} A \wedge T^{n_{1,0}} S^{n_{3,0}-n_{1,0}} A ) | n_{3,0} = n_{1,0} + n_{2,0} )$$
which simplifies (using the shift invariance) to $\P(A) = \eta$.  Thus we have $\eta < 3\eta/5$, a contradiction.  This proves Theorem \ref{fk-thm}.
\endprf

\begin{remark}  At present, the hypergraph regularity method is known to yield the Furstenberg-Katznelson recurrence theorem, but more powerful recurrence theorems, such as the Bergelson-Leibman polynomial recurrence theorem, the
Furstenberg-Katznelson IP-Szemer\'edi theorem, and the Furstenberg-Katznelson density Hales-Jewett theorem, have not yet been successfully obtained by this method (either in the finitary or infinitary settings).  It is not clear to the author whether this represents any fundamental limitations to the method.  A possible test problem would be the refinement of Szemer\'edi's theorem that the set of possible differences amongst the arithmetic progressions of a given length is syndetic (has bounded gaps); this was established for instance in \cite{furst} by ergodic methods but does not currently have a non-ergodic proof.
\end{remark}

\end{document}